\documentclass{article}
\usepackage{amsmath, amsxtra,amssymb, latexsym, amsfonts} 
\usepackage{multicol,color}
\usepackage{bbm}
\usepackage{float}
\usepackage{soul}
\usepackage{graphicx}
\definecolor{labelkey}{rgb}{0,0.08,0.45}
\definecolor{refkey}{rgb}{0,0.6,0.0}
\definecolor{Brown}{rgb}{0.45,0.0,0.05}
\definecolor{lime}{rgb}{0.00,0.8,0.0}
\definecolor{lblue}{rgb}{0.5,0.5,0.99}
\usepackage{pstricks-add}
\usepackage{amssymb}
\usepackage{theorem}
\usepackage{fancyhdr}
\usepackage{tikz}
\usepackage{enumerate}
\usepackage[margin=1.1in]{geometry}
\usetikzlibrary{arrows}
\usepackage{hyperref}
\usepackage{esint }
\definecolor{labelkey}{rgb}{0.6,0.6,0.6}
\definecolor{refkey}{rgb}{0,0.6,0.0}
\def\disp{\displaystyle}
\def\e{\epsilon}

\def\vt{\vartheta}

\def\O{\Omega}

\def\({\left(}
\def\){\right)}
\def\[{\left[}
\def\]{\right]}
\def\n{\left \|}
\def\en{\right \|}
\def\nn{\left \{ }
\def\hnn{\right \}}

\def\ox{\bar{x}}

\def\ou{\bar{u}}

\def\oT{\overline{T}}

\def\gg{\gamma}

\def\hat{\widehat}

\def\la{\langle}
\def\ra{\rangle}

\def\B{I\!\!B}

\def\R{\mathbb{R}}
\def\N{\mathbb{N}}

\def\co{\mbox{\rm co}\,}

\def\O{\Omega}

\def\vph{\varphi}
\def\emp{\emptyset}

\def\gg{\gamma}

\def\be{\beta}

\def\ph{\varphi}

\def\N{I\!\!N}
\def\th{\theta}

\newtheorem{theorem}{Theorem}[section]

\newtheorem{definition}[theorem]{Definition}
\theoremstyle{plain}{\theorembodyfont{\rmfamily}
}
\theoremstyle{plain}{\theorembodyfont{\rmfamily}
}
\theoremstyle{plain}{\theorembodyfont{\rmfamily}
}
\theoremstyle{plain}{\theorembodyfont{\rmfamily}
}

\theoremstyle{plain}{\theorembodyfont{\rmfamily}
}

\def\eq{\begin{equation}}
\def\eeq{\end{equation}}
\begin{document}
\begin{center}
\end{center}
\begin{center}
{\bf OPTIMAL CONTROL OF SWEEPING PROCESSES\\ IN UNMANNED SURFACE VEHICLE AND NANOPARTICLE MODELING}\\[3ex]
BORIS S. MORDUKHOVICH\footnote{Department of Mathematics, Wayne State University, Detroit, Michigan 48202, USA (boris@math.wayne.edu). Research of this author was partially supported by the USA National Science Foundation under grants DMS-1007132 and DMS-1512846, by the Australian Research Council under Discovery Project DP-190100555, and by Project~111 of China under grant D21024.} 
\quad DAO NGUYEN\footnote{Department of Mathematics and Statistics, San Diego State University, San Diego, CA 92182, USA (dnguyen28@sdsu.edu). Research of this author was supported by the AMS-Simon Foundation.}
\quad TRANG NGUYEN\footnote{Department of Mathematics, Wayne State University, Detroit, Michigan 48202, USA (daitrang.nguyen@wayne.edu). Research of this author was partially supported by the USA National Science Foundation under grant DMS-1512846.}
\end{center}
\begin{center}
{\bf Dedicated to Lionel Thibault, with deep respect and admiration} 
\end{center}
\small{\bf Abstract.} This paper addresses novel applications to practical modeling of the newly developed theory of necessary optimality conditions in controlled sweeping/Moreau processes with free time and pointwise control and state  constraints. Problems of this type appear, in particular, in dynamical models dealing with unmanned surface vehicles (USVs) and nanoparticles. We formulate optimal control problems for a general class of such dynamical systems and show that the developed necessary optimality conditions for constrained free-time controlled sweeping processes lead us to designing efficient procedures to solve practical models of this class. Moreover, the paper contains numerical calculations of optimal solutions to marine USVs and nanoparticle models in specific situations. Overall, this study contributes to the advancement of optimal control theory for constrained sweeping processes and its practical applications in the fields of marine USVs and nanoparticle modeling.\\[1ex]
{\em Key words.} Sweeping processes, optimal control, marine surface vehicles, nanoparticles, variational analysis, necessary optimality conditions,  discrete approximations, generalized differentiation\\[1ex]
{\em AMS Subject Classifications.} 49J52, 49J53, 49K24, 49M25, 90C30\vspace{-0.15in}

\section{Introduction}\label{intro}
\setcounter{equation}{0}\vspace*{-0.1in}

Sweeping process models, initially termed ``processus de rafle'', were pioneered by Jean-Jacques Moreau in the 1970s as a means to elucidate the dynamics inherent in elastoplasticity and related mechanical domains; see, e.g. \cite{moreau}. Such models, described in the form of the dissipative differential inclusions of the type 
\begin{equation}\label{sweep}
\left\{\begin{matrix}
\dot{x}(t)\in-N\big(x(t);C(t)\big)\;\textrm{ a.e. }\;t\in[0,T],\\x(0)=x_0\in C(0)\subset\R^n,
\end{matrix}\right.
\end{equation}
via the normal cone to the convex moving sets $C(t)$, have been expounded upon in many publications by Moreau and other researchers. We list here just a few of them \cite{aht,mb,brog,BrS,CT,hb1,KP,Kr,maury, venel,MM,smb,thi,tolstonogov}, where the reader can find further references to theoretical developments and various applications to aerospace, process control, robotics, bioengineering, chemistry, biology, economics, finance, management science, and engineering.

Lionel Thibault, a colleague of Moreau, has made fundamental contributions to the theory of sweeping processes in finite and infinite dimensions; see, e.g., \cite{aht,et,thi,thibault-book} and the references therein.

Since the Cauchy problem for sweeping processes of type \eqref{sweep} admits a {\em unique} solution under natural assumptions, there is no room for optimization in the framework of \eqref{sweep}. More recently, some {\em optimization problems} have been formulated and investigated for various {\em controlled sweeping processes}, where control actions enter moving sets, adjacent differential equations, and/or perturbations of the sweeping dynamics. We refer the reader to the stream of rather recent publications
\cite{ao,ac,bk,b1,CaoMordukhovich2018DCDS,chhm1,chhm,cmn18b,
cmn18a,pfs,rosario,hjm,michele,vera,zeidan} and the bibliographies therein, where {\em necessary optimality conditions} have been established and applied to models of friction and plasticity, robotics, traffic equilibria, ferromagnetism, hysteresis, and other fields of applied sciences. 

It should be mentioned that deriving necessary optimality conditions for controlled sweeping problems is much more challenging in comparison with optimal control of ODEs and Lipschitzian differential inclusions of the type $\dot x\in F(x)$; see, e.g., \cite{cl,Mord1,v}. The main difficulties come from the intrinsic {\em discontinuity} of the sweeping dynamics described by the normal cone mappings and the intrinsic presence of {\em pointwise state} as well as {\em irregular mixed} constraints on control and state functions.

This paper is devoted to applications of the novel theory of necessary optimality conditions developed in \cite{cmnn23} for a new class of controlled sweeping processes with free time and hard control-state constraints. For the first time in the literature, we address here practical optimal control models, which naturally arise in the dynamical modeling of marine unmanned surface vehicles and nanoparticles. In fact, the theoretical developments in \cite{cmnn23} have been largely motivated by the applications presented in this paper. 

In Section~\ref{problem}, we formulate and discuss the free-time sweeping optimal control problem of our study with overviewing some major results on discrete approximations and necessary optimality conditions distilled from \cite{cmnn23}. Employing these developments, we consider in Section~\ref{model} a broad class of planar dynamical models, which can be delineated via the free-time sweeping optimal control problem investigated in Section~\ref{problem}. Sections~\ref{marine} and \ref{sec:example2} formulate and study in detail optimal control problems for marine USVs and nanoparticles, respectively, with calculating optimal solutions to these models by using the obtained necessary optimality conditions for controlled sweeping constrained systems. The concluding Section~\ref{sec:Conclusions} summarizes the underlying features of the developed models and obtained results with discussing some directions of the future research.\vspace*{-0.15in}

\section{Sweeping Optimal Control with Free Time}\label{problem}

First we formulate a general free-time optimal control problem for the sweeping dynamics and highlight (with additional discussions and comments) some pivotal findings from \cite{cmnn23} that are essential for the subsequent applications to practical modeling. By $(P)$, we label the {\em sweeping optimal control problem} described as follows:
\begin{equation*}
\text{minimize}\quad J[x,u,T]:=\varphi(x(T),T)
\end{equation*}
over measurable controls $u(\cdot)$ and absolutely continuous sweeping trajectories $x(\cdot)$ defined on the variable time interval $[0,T]$ and satisfying the dynamic, control, and endpoint constraints 
\begin{equation}\label{Problem}
\left\{\begin{matrix}
\dot{x}(t)\in-N\big(x(t);C(t)\big)+g\big(x(t),u(t)\big)\;\textrm{ a.e. }\;t\in[0,T],\;x(0)=x_0\in C(0)\subset\R^n,\\
u(t)\in U\subset\R^d\;\textrm{ a.e. }\;t\in[0,T],\\ (x(T), T)\in \Omega_x\times \Omega_T\subset\R^n\times [0, \infty ), 
\end{matrix}\right.
\end{equation}
where $\O_x$ and $\O_T$ are subsets of $\R^n$ and $[0,\infty)$, respectively, and where $C(t)$ is a convex polyhedron given by
\begin{equation}\label{C}
\left\{\begin{matrix}
C(t):=\bigcap_{j=1}^{s}C^j(t)\textrm{ with }C^j(t):=\nn x\in\R^n\;\big|\;\la x^j_*(t),x\ra\le c_j(t)\hnn,
\\ 
\|x_\ast^j(t)\|=1$, $j=1,\ldots,s,\,t\in [0,\infty).
\end{matrix}\right.
\end{equation}
Recall that the normal cone to the convex set in \eqref{Problem} is defined by
\begin{equation*}
N(x;C):=\big\{v\in\R^n\;\big|\;\la v,y-x\ra\le 0,\;y\in C\big\}\textrm{ if }x\in C\textrm{ and }N(x;C):=\emp\textrm{ if }x\notin C.
\end{equation*}
The latter tells us that problem $(P)$ automatically contains the {\em pointwise state constraints}
\begin{equation*}
x(t)\in C(t),\textrm{ i.e., }\la x^j_*(t),x(t)\ra\le c_j(t)\;\textrm{ for all }\;t\in [0,T]\;\textrm{ (with different } T)\;\textrm{ and }\;j=1,\ldots,s.
\end{equation*}
In fact, the sweeping dynamics intrinsically induce {\em irregular mixed constraints} on controls and trajectories that are the most challenging and largely underinvestigated even in classical control theory for smooth ODEs.

All the triples $(x(\cdot),u(\cdot),T)\in W^{1,2}([0,T];\R^n)\times L^2([0,T];\R^d)\times\R_+$ satisfying \eqref{Problem} are called {\em feasible solutions} to $(P)$. We identify the trajectory $x:[0,T]\to\R^n$ with its extension to the interval $(0,\infty)$ given by
\begin{equation*}
x_e(t):=x(T)\;\textrm{ for all }\;t>T,
\end{equation*}
and for $x(\cdot)\in W^{1,2}([0,T],\R^n)$ define the norm 
\begin{equation*}
\n x\en_{W^{1,2}}:=\n x(0)\en+\n \dot{x}_e\en_{L_2}.
\end{equation*}

Let us now specify what we mean by a {\em $W^{1,2}\times L^2\times \R_+$-local minimizer} for (P) and its relaxation; cf.\ \cite{cmn18a} when the duration of the process is fixed.

\begin{definition}\label{Def3.1}
A feasible solution $(\ox(\cdot),\ou(\cdot),\oT)$ to problem $(P)$ qualifies as a {\sc $W^{1,2}\times L^2\times \R_+$}-local minimizer for this problem if there exists $\e>0$ such that $J[\ox,\ou,\oT]\le J[x,u,T]$ holds for all feasible solutions $(x(\cdot),u(\cdot),T)$ to $(P)$ satisfying the localization condition
\begin{equation*}
\int_0^{\oT}\(\n\dot{x}_e(t)-\dot{\ox}_e(t)\en^2+\n u(t)-\ou(t)\en^2\)dt+(\oT-T)^2<\e.
\end{equation*}
\end{definition}

The {\em relaxed} version $(R)$ of problem $(P)$ is constructed as follows:
\begin{equation*}
\textrm{ minimize }\hat{J}[x,u,T]:=\varphi(x(T),T)
\end{equation*}
over absolutely continuous trajectories of the convexified differential inclusion
\begin{equation}\label{conv}
\dot{x}(t)\in-N\big(x(t);C(t)\big)+\co g\big(x(t),U\big)\;\textrm{ a.e. }\;t\in[0,T],\;x(0)=x_0\in C(t)\subset\R^n,
\end{equation}
where ``co" represents the convex hull of the set. 

\begin{definition}\label{relaxed} Let $(\ox(\cdot),\ou(\cdot),\oT)$ be a feasible solutions for problem $(P)$. We say that it is a
{\sc relaxed $W^{1,2}\times L^2\times\R_+$-local minimizer} for $(P)$ if the following condition holds: there exists $\e>0$ such that
\begin{equation*}
\ph\big(\ox(\cdot),\oT\big)\le\ph\big(x(\cdot),T\big)\;\textrm{ whenever }\;\int_0^{\oT}\(\n\dot{x}_e(t)-\dot{\ox}_e(t)\en^2+\n u(t)-\ou(t)\en^2\)dt+(\oT-T)^2<\e,
\end{equation*}
where $u(t)\in\co U$ a.e.\ on $[0,T]$ with $u(\cdot)$ representing a measurable control, and $x(\cdot)$ denotes a relaxed trajectory of the convexified inclusion \eqref{conv}, which can be uniformly approximated in $W^{1,2}([0,T];\R^n)$, by feasible trajectories to $(P)$ generated by piecewise constant controls $u^k(\cdot)$ on $[0,T]$, the convex combinations of which converge strongly to $u(\cdot)$ in the norm topology $L^2([0,T];\R^d)$.
\end{definition}\vspace*{-0.05in}

It follows from the construction of the relaxed problem $(R)$ that the local minimizers in Definitions~\ref{Def3.1} and \ref{relaxed} agree under the convexity assumptions on the data of $(P)$. Due to the nonatomicity of the Lebesgue measure, such a relaxation stability phenomenon also holds in broad nonconvex settings; see \cite{cmn18a,cmnn23} for more details. 

Now we formulate the {\em standing assumptions} on the given data of the sweeping control system in \eqref{Problem} and \eqref{C}, which are all satisfied in the practical models studied in the subsequent sections.\\[1ex]
{\bf(H1)} The control set $U\ne\emp$ is closed and bounded in $\R^d$.\\[1ex]
{\bf(H2)} The generating functions $x^j_{\ast}(\cdot)$ and $c_j(\cdot)$ of the moving polyhedron in \eqref{C} are Lipschitz continuous with a common Lipschitz constant $L$.\\[1ex]
{\bf(H3)} The {\em uniform Slater condition} is fulfilled:
\begin{equation}\label{sla}
\mbox{for every }\;t\in[0,T]\;\mbox{  there exists }\,x\in\R^n\;\mbox{ such that }\;\la x^j_{\ast}(t), x\ra<c_j(t)\;\mbox{ whenever }\;i=1,\ldots,s.
\end{equation}
{\bf(H4)} The perturbation mapping $g\colon\R^n\times U\to\R^n$ exhibits {\em Lipschitz continuity} with respect to $x$ uniformly on $U$ whenever $x$ belongs to a bounded subset of $\R^n$, and the {\em sublinear growth condition} holds: there is $\be\ge 0$ with
\begin{equation*}
\|g(x,u)\|\le\be\big(1+\|x\|\big)\;\mbox{ for all }\;u\in U.
\end{equation*}
{\bf(H5)} The endpoint-final time constraint set $\O_x\times\O_T$ is closed.

Note that the uniform Slater condition \eqref{sla}, introduced recently in \cite{hjm} always holds under the fulfillment of the {\em linear independence constraint qualification} (LICQ)
\begin{equation}\label{plicq}
\Big[\sum_{j\in I(t,x)}\lambda_jx^j_*(t)=0,\;\lambda_j\in\R\Big]\Longrightarrow\big[\lambda_j=0\;\textrm{ for all }\;j\in I(t,x)\big],
\end{equation}
where the moving active index set at $(t,x)$ is defined by
\begin{equation*}
I(t,x):=\big\{j\in\{1,\ldots,s\}\;\big|\;\la x^j_{\ast}(t),x\ra=c_j(t)\big\}.
\end{equation*}
On the other hand, \eqref{sla} yields the {\em positive linear independence constraint qualification} (PLICQ) corresponding to \eqref{plicq} with $\lambda_j\in\R$ replaced by $\lambda_j\in\R_+$ for $j\in I(t,x)$. It follows from \cite[Theorem~2]{hjm} that the assumptions imposed above ensure that any feasible control $u(\cdot)$ to $(P)$ generates the unique Lipschitz continuous trajectory $x(\cdot)$ of the sweeping process in \eqref{Problem}. 

Our approach to the study and solving dynamic optimization problems is based on the {\em method of discrete approximations} developed by the first author \cite{m95,Mord1} for Lipschitzian differential inclusions. Since the sweeping dynamics is highly non-Lipschitzian, applications of this method to controlled sweeping processes require significant modifications, which have been done in \cite{CCMN21,cm2,b1,cg21,chhm,cmn18a,cmnn23,mn20} for various classes of sweeping control problems. To proceed in the setting of \eqref{Problem} and \eqref{C} of $(P)$, for each natural number $k\in\N:=
\{1,2,\ldots\}$ consider the time $T_k$ approximating $T$ and form the uniform grid on $[0,T]$ by
\begin{equation}\label{grid}
 t^k_0=0,\; t^k_{i+1}=t^k_i+h^k\;\mbox{ as }\;i=0,1,\ldots,k-1,\;\;t^k_k:=T_k,\;h^k:=k^{-1}T_k.
\end{equation}

The following result, distilled from \cite[Theorem~2.3]{cmnn23}, verifies a {\em strong approximation} of {\em any feasible solution} to the continuous-time sweeping control problem $(P)$ by a sequence of feasible solutions to perturbed {\em discrete-time sweeping processes}. Beside being useful in deriving necessary optimality conditions for $(P)$, this result justifies the possibility to replace the infinite-dimensional sweeping dynamics by its finite-dimensional counterparts in {\em numerical calculations}, which is important for practical modeling. \vspace*{-0.05in}

\begin{theorem}\label{Thr1}
Let $(\ox(\cdot),\ou(\cdot),\oT)\in W^{1,2}([0,\oT];\R^n)\times L^2([0,\oT];\R^d)$ be an arbitrary feasible solution to problem $(P)$ under the imposed standing assumptions. Then we have the assertions:\\\vspace*{0.05in}
{\bf(i)} There exists a sequence of piecewise constant functions $\{u^k(\cdot)\;|\;k\in \N \}$ such that $u^k(\cdot)\to\ou(\cdot)$ strongly
in the $L^2$-norm topology on $[0,\oT]$ as $k\to\infty$.\\\vspace*{0.05in}
{\bf(ii)} There exists a sequence of piecewise linear functions $\{x^k(\cdot)\:|\; k\in \N \}$, which converges strongly to 
$\ox(\cdot)$
in the $W^{1,2}$-norm topology on $[0,\oT]$ with $x^k(0)=\ox(0)$
for all $k\in\N$ and such that
\begin{equation*}
\dot x^k(t)\in -\mbox{N}(x^k(t^k_i);C^k_i)+g(x^k(t^k_i), u^k(t^k_i))+\tau^k_i\B,\;t\in[t^k_i,t^k_{i+1}),
\end{equation*}
for $i=0,\ldots,k-1$, where $\tau^k_i\geq 0$, $h^k\sum^{k-1}_{i=1}\tau^k_i\to 0$ as $k\to \infty$, $\B$ stands for the closed unit ball, and the perturbed polyhedra $C^k_i$ are defined by
\begin{equation*}
C^k_i:=\disp\bigcap_{j=1}^s C^k_{ij}\;\mbox{ with }\;C^k_{ij}:=\big\{x\in\R^n\;\big|\;\la x,x^j_{\ast}(t^k_i)\ra \le c^{ij}_k\big\}
\end{equation*}
with some vectors $x^j_{\ast}(t^k_i)$ and numbers $c^{ij}_k$. Furthermore, all the approximating arcs $x^k(\cdot)$ are Lipschitz continuous on $[0,\oT]$ with the same Lipschitz constant as $\bar{x}(\cdot)$.\\\vspace*{0.05in}
{\bf(iii)} The piecewise linear extensions of $x^j_{\ast k}(t^k_j)$ and $c^{ij}_k$ on the continuous time interval $[0,\oT]$ converge uniformly on $[0,\oT]$ as $k\to \infty$ to $c_j(t)$ and $x^j_{\ast}(t)$, respectively. 
\end{theorem}

Having in hand the strong approximation results of Theorem~\ref{Thr1} applied to the designated relaxed $W^{1,2}\times L^2\times\R_+$-local minimizer $(\ox(\cdot),\ou(\cdot),\oT)$ of $(P)$, we construct the sequence of {\em discrete approximation problems} $(P_k)$ with a varying grid whose {\em optimal solutions exists} (assuming that the cost function $\ph$ is lower semicontinuous) and {\em strongly converge} in $W^{1,2}\times L^2\times\R_+$-norm to
$\{\bar{x}(\cdot),\bar{u}(\cdot),\oT\}$. For each $k\in\N$, define problem $(P_k)$ as follows:
\begin{eqnarray*}
\mbox{minimize }\;J_k[x^k,u^k,T_k]:=\varphi(x^k_k,T_k)+(T_k-\oT)^2+
\end{eqnarray*}
$$\sum_{i=0}^{k-1}\int_{t^k_i}^{t^k_{i+1}}\left(\n \frac{x^k_{i+1}-x^k_i}{h^k}-\dot{\bar{x}}(t)\en^2+\n u^k_i-\bar{u}(t)\en^2\right)dt$$
over $(x^k,u^k,T_k):=(x^k_0,x^k_1,\ldots,x^k_{k-1}, x^k_k,u^k_0,u^k_1,\ldots,u^k_{k-1},T_k)$ satisfying the constraints
\begin{equation}\label{gph}
x^k_{i+1}- x^k_i\in -h^kF^k_i(t^k_i,x^k_i,u^k_i)\;
\textrm{ for }\;i=0,\ldots,k-1,
\end{equation}
\begin{equation*}
x^k_0:=x_0\in C(0),
\end{equation*}
\begin{equation*}
(x^k_k, T_k)\in \O^k_x\times \O^k_T:= (\O_x+\delta^k\B)\times (\O_T+\delta^k),
\end{equation*}
\begin{equation*}
\sum_{i=0}^{k-1}\int_{t^k_i}^{t^k_{i+1}}\(\n\frac{x^k_{i+1}-x^k_i}{h^k}-\dot{\ox}(t)\en^2+\n u^k_i-\ou(t)\en^2\)dt\le \e,
\end{equation*}
\begin{equation*}
\left\| \frac{x^k_{i+1}-x^k_i}{h^k}\right\| \le L +1\;\text{ for all } i=0,\ldots ,k-1,
\end{equation*}
\begin{equation*}
u^k_i\in U\;\textrm{ for }\;i=0,\ldots,k-1,
\end{equation*}
\begin{equation*}
\|((x^k_i,u^k_i)-(\ox(t^k_i),\ou(t^k_i))\|\le \e \textrm{ for } i=0,\ldots,k-1,
\end{equation*}
\begin{equation*}
|T_k-\oT|\le \e
\end{equation*}
\begin{equation*}
\la x^j_\ast(t), x^k_k\ra \leq c_j(t)\textrm{ for all } j =1,\ldots,s,
\end{equation*}
where the discrete velocity mappings $F^k_i$ are given by
\begin{equation*}
F^k_i(t^k_i,x^k_i,u^k_i):=N(x^k_i;C^k_i)-g(x^k_i, u^k_i)-\tau^k_i\B,
\end{equation*}
where $\delta^k:=\|\ox(\oT)-\hat x^k(\oT)\|$ with the sequence 
$\{x^k(\oT)\}$ constructed in Theorem~\ref{Thr1} for $\ox(\cdot)$, where $\e >0$ is taken from Definition~\ref{relaxed} for $\ox(\cdot)$, and where $L$ is the Lipschitz constant of $\ox(\cdot)$ on $[0,\oT]$.

Employing Theorem~\ref{Thr1} and using the construction of problems $(P_k)$ allow us to verify the strong $W^{1,2}\times L^2\times\R_+$ convergence of the extended (as in Theorem~\ref{Thr1}) optimal solutions to the designated local minimizer $(\ox(\cdot,\ou(\cdot),\oT)$ of $(P)$. This means that optimal solutions to the discrete approximation problems $(P_k)$ can be viewed as {\em suboptimal solutions} to the original sweeping control problem $(P)$. 

To derive further precise {\em necessary optimality conditions} for the given local minimizer in $(P)$, we develop the following {\em two-step procedure}:\\[1ex]
{\bf Step~A:} {\em Obtain necessary conditions for optimal solutions to $(P_k)$}\\[1ex]
{\bf Step~B:} {\em By passing to the limit in the necessary optimality conditions for $(P_k)$ as $k\to\infty$, establish necessary optimality conditions for the designated local minimizers $(\ox(\cdot),\ou(\cdot),\oT)$ in $(P)$}.

The goal of Step~A is accomplished by reducing each problem $(P_k)$ to a mathematical program with geometric and functional constraints and deriving necessary optimality conditions in the latter setting by employing appropriate tools of {\em generalized differentiation} in {\em variational analysis}. Major requirements to such tools to be useful for furnishing the limiting procedure in Step~B include {\em robustness} with respect to small perturbations, comprehensive {\em calculus rules} of generalized differentiation, and the ability to deal with geometric constraints of the {\em graphical type} \eqref{gph}, which are crucial to reflect the sweeping dynamics. Note to this end that Clarke's constructions of generalized differentiation \cite{cl} are not suitable for these purposes, since his normal cone is too large for graphical sets associated with functions and mappings even in very simple situations as, e.g., for the graph of $f(x):=|x|$, where Clarke's normal cone at $(0,0)$ is the entire plane $\R^2$. On the other hand, the {\em limiting generalized differential constructions}, introduced by the first author and then developed in numerous publications (see, e.g., the books \cite{Mord,m-book2,rw,thibault-book} and the references therein), provide the desired variational machinery for passing to the limit in Step~B and deriving in this way an adequate collection of necessary optimality conditions in the sweeping control problem $(P)$ that are satisfactory for applications to the practical models considered below. 

The underlying feature of the sweeping dynamics and its discrete approximations is their descriptions via the normal cone mapping, which accumulates first-order information about the process. Thus the {\em adjoint systems} in both discrete-time and 
continuous-time problems unavoidably require a dual-type generalized derivative ({\em coderivative}) of normal cone mappings that are formalized by the first author as the {\em second-order subdifferential} (or {\em generalized Hessian}).
The {\em explicit calculations} of this second-order construction for set-valued mappings, which appear in the sweeping dynamics \eqref{Problem} and \eqref{gph}, play a crucial in the realization of the method of discrete approximations for optimal control of sweeping processes. The reader can find such calculations and related results in the aforementioned papers on sweeping optimal control with the references therein. A comprehensive theory of second-order subdifferentials with a wide spectrum of applications is presented in \cite{m-book2nd}.

To formulate the major necessary optimality conditions for problem $(P)$ proved in \cite[Theorem~5.2]{cmnn23} by using the method of discrete approximations, we need the following notion. Observe to this end that Motzkin's theorem of the alternative gives us the representation of the normal cone $N\(x(t);C(t)\)$ as $\{\sum_{j\in I(t,x)}\eta^j(t) x^j_*(t)\;|\;\eta^j(t)\ge 0\}$.\vspace*{-0.1in}

\begin{definition}\label{active cone}
Let $x(\cdot)$ be a solution to the controlled sweeping process \eqref{Problem}, i.e.,
\begin{equation*}
-\dot{x}(t)=\sum_{j=1}^s\eta^j(t)x^j_*(t)-g\big(x(t),u(t)\big)\;\textrm{ for all }\;t\in [0,T),
\end{equation*}
where $\eta_j\in L^2([0,T];\mathbb{R}^+)$ and $\eta_j(t) =0$ for a.e. $t$ such that $j\not\in I(t,x(t))$. We say that the normal cone to
$C(t)$ is active along $x(\cdot)$ on the set $E\subset[0,T]$ if for a.e. $t\in E$ and all $j\in I(t,x(t))$ it holds that
$\eta_j(t)>0$. Denoting by $E_0$ the largest subset of $[0,T]$ where the normal cone to $C(t)$ is active along $x(\cdot)$ 
$($i.e., the union of the density points of $E$ such that the normal cone to $C(t)$ is active with respect to this set$)$, we simply say that the normal cone is active along $x(\cdot)$ provided that $E_0=[0,T]$.
\end{definition}\vspace*{-0.1in}

Note that the requirement that the normal cone is active along a trajectory of \eqref{Problem} distinguishes \eqref{Problem} from trajectories of the classical controlled ODE $\dot{x}=g(x,u)$ which satisfy the {\em tangency condition} $\langle g(x(t),u(t)), x_\ast^j(t)\rangle \le 0$ for a.e. $t$ such that $j\in I(t,x(t))$ and thus $x(t)\in C(t)$ for all $t\in E$. Observe also that the above requirement holds automatically when the set $I(t,x(t))$ is empty, i.e., when $x(t)$ stays in the interior of $C(t)$.\vspace*{-0.1in}

\begin{theorem}\label{Thr7}
Let $(\ox(t),\ou(t),\oT),\;0\leq t\leq\oT$, be a relaxed $W^{1,2}\times L^2\times\R_+$--local minimizer to problem $(P)$. In addition to the standing assumptions, suppose that LICQ holds along $\ox(t)$, $t\in[0,T]$, and that $\ph$ is locally Lipschitzian around $(\bar{x}(\oT),\oT)$. Then there exist a multiplier $\mu\ge 0$, a nonnegative vector measure 
$\gg_>=(\gg_>^1,\ldots,\gg_>^s)\in C^*([0,\oT];\R^s)$, and a signed vector measure 
$\gg_0=(\gg_0^1,\ldots,\gg_0^s)\in C^*([0,\oT];\R^s)$ together with adjoint arcs $p(\cdot)\in W^{1,2}([0,\oT];\R^n)$ and $q(\cdot)\in BV([0,\oT];\R^n)$ such that the following conditions are fulfilled:\\[1ex]
$\bullet$ The sweeping trajectory representation
\begin{equation}\label{37}
-\dot{\ox}(t)=\sum_{j=1}^s\eta^j(t)x^j_*(t)-g\big(\ox(t),\ou(t)\big)\;\textrm{ for a.e. }\;t\in [0,\oT),
\end{equation}
where the functions $\eta^j(\cdot)\in L^2([0,\oT]);\R_+)$ are uniquely determined for a.e.\ $t\in[0,\oT)$ by representation \eqref{37}.\\[1ex]
$\bullet$ The adjoint arc inclusion
\begin{equation*}
\big(-\dot{p}(t), \psi(t)\big)\in \co\partial\la q(t), g\ra (\ox(t),\ou(t)\big)\;\textrm{ for a.e. }\;t\in[0,\oT],
\end{equation*}
where the limiting subdifferential is taken with respect to $(x,u)$, where $\psi(\cdot)\in L^2([0,\oT];\R^d)$ satisfies the 
\begin{equation*}
\psi (t)\in \co N(\ou(t);U)\;\textrm{ for a.e. }\;t\in[0,\oT],
\end{equation*}
where the right continuous representative of $q(\cdot)$ is
given by
\begin{equation*}
q(t)=p(t)-\int_{(t,\oT]} \sum_{j=1}^sd\gg^j(\tau)x^j_*(\tau)
\end{equation*}
for a.e. $t\in[0,\oT]$ except at most a countable subset, and where $\gg:=\gg_>+\gg_0$. Moreover, $p(\oT)=q(\oT)$.\\[1ex]
$\bullet$ The tangential maximization condition: if the limiting normal cone is generated as
\begin{equation*}
N(\ou(t);U)=T^*(\ou(t);U):=\big\{v\in\R^n\big|\;\la v,u\ra\le 0\;\mbox{ for all }\;u\in T\big(\ou(t);U\big)\big\}
\end{equation*}
by some tangent set $T(\ou(t);U)$ associated with $U$ at $\ou(t)$, then we have
\begin{equation*}
\big\la \psi(t),\ou(t)\big\ra=\max_{u\in T(\ou(t);U)}\big\la \psi(t),u\big\ra\;\textrm{ for a.e. }\;t\in[0,T].
\end{equation*}
In particular, the global maximization condition
\begin{equation*}
\big\la \psi(t),\ou(t)\big\ra=\max_{u\in U}\big\la \psi(t),u\big\ra\;\textrm{ for a.e. }\;t\in[0,\oT]
\end{equation*}
is satisfied provided that the control set $U$ is convex.\\[1ex]
$\bullet$ The dynamic complementary slackness conditions
\begin{equation*}
\big\la x^{j}_\ast(t),\ox(t)\big\ra<c_j(t)\Longrightarrow\eta^j(t)=0\;\mbox{ and }\;\eta^j(t)>0\Longrightarrow\;\big\la x^{j}_\ast(t),q(t)\big\ra=0
\end{equation*}
for a.e.\ $t\in[0,\oT]$ and all indices $j=1,\ldots,s$.\\[1ex]
$\bullet$ The transversality conditions at the optimal final time: there exist numbers $\eta^j(\oT)\ge 0$ whenever $j\in I(\oT,\ox(\oT))$ ensuring the relationships
\begin{equation}\label{42}
\begin{array}{ll}
\big(-p(\oT)-\disp\sum_{j\in I(\oT,\ox(\oT))}\eta^j(\oT) x^{j}_\ast(\oT),\bar{H}\big)\in\mu \partial\varphi(\bar{x}(\oT),\oT)+N((\ox(\oT),\oT);\O_x\times \O_T),\\
\eta^j(\oT)>0\Longrightarrow j\in I\big(\oT,\ox(\oT)\big),
\end{array}
\end{equation}
where $\bar{H}:={\oT}^{-1}\int_{0}^{\oT}\langle p(t),\dot{\bar{x}}(t)\rangle dt$ is a characteristic of the optimal time.\\[1ex]
$\bullet$ The endpoint complementary slackness conditions
\begin{equation*}
\big\la x^{j}_\ast(\oT),\ox(\oT)\big\ra<c_j(\oT)\Longrightarrow\eta^j(\oT)=0\;\mbox{ for all }\;j\in I\big(\oT,\ox(\oT)\big)
\end{equation*}
with the nonnegative numbers $\eta^j(\oT)$ taken from \eqref{42}.\\[1ex]
$\bullet$ The nonatomicity condition: If $t\in[0,\oT)$ and $\la x^{j}_\ast(t),\ox(t)\ra<c_j$ for all $j=1,\ldots,s$, then there
exists a neighborhood $V_t$ of $t$ in $[0,\oT)$ such that $\gg^j_0(V)=\gg^j_>(V)=0$ for all Borel subsets $V$ of $V_t$. In particular, $\mathrm{supp}(\gamma^j_>)$ and $\mathrm{supp}(\gamma^j_0)$ are contained in the set $\{t\;|\;j\in I(t,\ox(t))\}$.\\[1ex]
$\bullet$ The general nontriviality condition
\begin{equation*}
(\mu,p,\|\gamma_0\|_{TV},\|\gg_>\|_{TV})\ne 0,
\end{equation*}
accompanied by the support condition
\begin{equation*}
\mathrm{supp}(\gamma_>)\cap \mathrm{int}(E_0)=\emptyset,
\end{equation*}
which holds provided the normal cone is active on a set with nonempty interior.\\[1ex]
$\bullet$ The enhanced nontriviality: we have $\mu=1$ for the cost function multiplier in \eqref{42} provided that $\la x^j_\ast(t),\ox(t)\ra<c_j(t)$ for all $t\in[0,\oT]$ and all indices $j=1,\ldots,s$.
\end{theorem}\vspace*{-0.1in}

In the subsequent sections, we develop applications of the obtained results to practical models formulated in the form of the sweeping optimal control problem $(P)$.\vspace*{-0.1in}

\section{Modeling}\label{model}

In this section, we delineate a class of rather general dynamical optimization models on the plane whose dynamics are described by sweeping processes. These models deal with $n\ge 2$ objects that have arbitrary shapes identified as virtual safety disks of different radius $R_i$, $i=1,\ldots,n$, on the plane. Each object aims at reaching the target by the shortest path with the minimum time $\oT$ while avoiding the other $n-1$ static and/or dynamic obstacles. We define the configuration space of objects at some time $t\in[0,T]$ by the vector $x(t)=(x^1(t),\ldots,x^n(t))\in\R^{2n}$ with the variable ending time $T$, where $x^i(t)=(\|x^i(t)\|\cos\th^x_i,\|x^i(t)\|\sin\th^x_i) \in\R^2$ denotes the Cartesian position of the $i$-th object, and where $\th^x_i$ stands for the corresponding constant direction which is the smallest positive angle in standard position formed by the positive $x$-axis and vectors $x^i$ with $0\in\R^2$ as the target. The configuration is admissible when the motion of different objects is safe by imposing the 
{\em noncollision/nonoverlapping condition}. This can be formulated mathematically as
\begin{eqnarray*}
A:=\big\{x=\(x^1,\ldots,x^n\)\in\mathbb{R}^{2n}\big|\;D_{ij}(x)\ge 0\;
\mbox{ for all }\;i,j\in\{1,\ldots,n\}\big\},
\end{eqnarray*}
where $D_{ij}(x):=\|x^{i}-x^j\|-(R_i+R_j)$ is the distance between the disks $i$ and $j$. 

Describing the motion of objects without collisions can be outlined as follows: starting from an admissible configuration at time $t_k\in[0,T]$ (with different $T$), consider $x_k:=x(t_k)\in A$. Then the next configuration after the period of time $h>0$ is $x_{k+1}=x(t_k+h)$. To ensure a nonoverlapping motion of all objects at the next time $t_k+h$ for a small value of $h>0$, the next configuration should also be admissible, i.e., $x(t_k+h)\in A$. This implies that the constraint $D_{ij}(x(t_k+h))\ge 0$ should be satisfied. To verify the latter, by using the first-order Taylor expansion at $x_k\ne 0$ we deduce the constraint on the velocity vector as
\begin{equation}\label{Dij}
D_{ij}\big(x(t_k+h)\big)=D_{ij}\big(x(t_k)\big)+h\nabla D_{ij}\big(x(t_k)\big)\dot{x}(t_k)+\;o(h),\;\textrm{ for small}\;h>0.
\end{equation}
where $\nabla D_{ij}$ is the gradient of the distance function. This constraint will be used to construct the next configuration or the next reference position of objects in order to avoid static and/or dynamic obstacles. In this regard, set the vector $V(x)$ to be the desired velocity of all objects.
The admissible velocities providing no collision during the navigation of objects are defined as
\begin{eqnarray*}
C_h(x):=\big\{V(x)\in\R^{2n}\;\big|\;D_{ij}(x)+h\nabla D_{ij}(x)V(x)\ge 0\;\mbox{ for all }i,j\in\{1,\ldots,n\},\;i<j\big\},\quad x\in\R^{2n}.
\end{eqnarray*}
Taking now the admissible velocity $\dot{x}(t_k)\in C_h(x_k)$ gives us
\begin{eqnarray*}
D_{ij}(x_k)+h\big\la\nabla D_{ij}(x_k),\dot{x}(t_k)\big\ra\geq 0.
\end{eqnarray*}
Skipping the term $o(h)$ for small $h$, it follows from \eqref{Dij} that $D_{ij}(x(t_k+h))\geq 0$, i.e., $x(t_k+h)\in A$.
Since all the objects intend to reach the target by the {\em shortest path}, by taking into account that in the absence of obstacles the objects tend to keep their spontaneous velocities till reaching the target and that $\|x^{i}(t)-x^j(t)\|=R_i+R_j$ if the designated object
touches the obstacles, we conclude that the object's velocity should be adjusted in order to keep the distance to be at least $R_i+R_j$ by using some {\em control actions} in the velocity term. The desired {\em controlled} velocities can be now represented as
\begin{eqnarray*}
g\big(x(t),u(t)\big)=\big(s_1\|u^1(t)\|\cos\th^u_1(t),s_1\|u^1(t)\|\sin\th^u_1(t),
\ldots,s_n\|u^n(t)\|\cos\th^u_n(t),s_n\|u^n(t)\|\sin\th^u_n(t)\big),
\end{eqnarray*}
where $\th^u_i$ stands for the corresponding constant direction, which is the smallest positive angle in standard position formed by the positive $x$-axis and vectors $u^i(t)$, and where $s_i\leq 0$ denotes the speed of the object $i$, with practically motivated control constraints represented by
\begin{equation}\label{ut}
u(t)=\big(u^1(t),\ldots,u^n(t)\big)\in U\;\mbox{ for a.e. }t\in[0,T]
\end{equation}
with the control set $U\subset\R^n$ to be specified below in particular settings.

Having the above discussions in mind, let us describe the model dynamics as a sweeping process. To proceed, for any $k\in \N$ consider the ending time $T_k$ and the grid as in \eqref{grid} with $x_i^k:=x^k(t^k_i)$ for $i=1,\ldots,k$. According to the model description, we have the algorithm
\begin{eqnarray*}
x^k_0\in A\;\mbox{ and }\;x^k_{i+1}:=x^k_i+h^kV^k_{i+1}\;\mbox{ for all }\;i=0,\ldots,k-1,
\end{eqnarray*}
where $V^k_{i+1}$ is defined as the {\em projection} of $g(x^k_i,u^k_i)$ onto the admissible velocity set $C_{h^k}(x^k_i)$ by
\begin{equation}\label{V_k}
V^k_{i+1}:=\disp\Pi\big(g(x^k_i,u^k_i);C_{h^k}(x^k_i)\big),\quad i=0,\ldots,k-1.
\end{equation}
Using the construction of $x^k_i$ for $0\le i\le k-1$ and $k\in\N$, design next a sequence of piecewise linear mappings $x^k\colon[0,T_k]\to\R^{2n}$, $k\in\N$, which pass through those points by
\begin{equation}\label{x_k}
x^k(t):=x^k_i+(t-t^k_i)V^k_{i+1}\;\mbox{ for all }\;t\in I^k_i:=[t^k_i,t^k_{i+1})
\end{equation}
with $i=0,\ldots,k-1$. Whenever $k\in\N$, we clearly have the relationships
\begin{equation}\label{xtki}
x^k(t^k_i)=x^k_i=\underset{t\to t^k_i}{\lim}x^k_i(t)\;\mbox{ and }\;\dot{x}^k_i(t):=V^k_{i+1}
\mbox{ for all } t\in(t^k_i,t^k_{i+1}).
\end{equation}
As discussed in \cite{mb}, the solutions to \eqref{x_k} in the {\em uncontrolled} setting of \eqref{V_k} with $g=g(x(t))$ uniformly converge on $[0,\oT_k]$ to a trajectory of a certain perturbed sweeping process. The {\em controlled} model under consideration here is significantly more involved. For all $x(t)\in\R^{2n}$, consider the set
\begin{equation}\label{Kx}
K(x(t)):=\big\{y(t)\in\R^{2n}\big|\;D_{ij}(x(t))+\nabla D_{ij}(x(t))(y(t)-x(t))\ge 0\;\mbox{ whenever }\;i<j\big\},
\end{equation}
which allows us to represent the algorithm in \eqref{V_k}, \eqref{x_k} as
\begin{eqnarray*}
x^k_{i+1}=\Pi\big(x^k_i+h^kg(x^k_i,u^k_i);K(x^k_i)\big)\;\mbox{ for }\;i=0,\ldots,k-1.
\end{eqnarray*}
Thus it can be equivalently rewritten in the form of
\begin{eqnarray*}
x^k\big(\vartheta^k(t)\big)=\Pi\big(x^k(\tau^k(t))+h^kg\big(x^k(\tau^k(t)),u^k(\tau^k(t)\big);K(x^k(\tau^k(t))\big)
\end{eqnarray*}
for all $t\in[0,\oT_k]$, where the functions $\tau^k(t):=t^k_i$ and $\vartheta^k(t):=t^k_{i+1}$ for all $t\in I^k_i$. Using the construction of the convex set $K(x)$ in \eqref{Kx} and the definition of the normal cone together with the relationships in \eqref{xtki}, we arrive at the sweeping process inclusions
\begin{equation}\label{dotxk}
\dot{x}^{k}(t)\in N\big(x^{k}(\vt^{k}(t));K(x^{k}(\tau^{k}(t)))\big)+g\big(x^{k}(\tau^{k}(t)),u^{k}(\tau^{k}(t))\big)\;\mbox{ a.e. }\;t\in[0,\oT]
\end{equation}
with $x^{k}(0)=x_0\in K(x_0)=A$ and $x^{k}(\vt^{k}(t))\in K(x^{k}(\tau^{k}(t)))$ on $[0,\oT]$. To formalize \eqref{dotxk} as a controlled perturbed sweeping process of type $\eqref{Problem}$, we define the convex polyhedron $C(t)\subset\R^{2n}$ by
\begin{equation}\label{setC}
C(t):=\bigcap\big\{x(t)\in\R^{2n}\;\big|\;\la x^j_*(t),x(t)\ra\le c_j(t),\;j=1,\ldots,n-1\big\}
\end{equation}
with $c_j(t):=-(R_{obj}+R_{obs})$, where $R_{obj}$, $R_{obs}$ are the radii of the designated object and the obstacle, respectively. The $n-1$ vertices of the polyhedron are given by
\begin{equation}\label{e}
x^j_*:=e_{j1}+e_{j2}-e_{(j+1)1}-e_{(j+1)2}\in\R^{2n}, \; j=1,\ldots,n-1,
\end{equation}
where $e_{ji}$ ($j=1,\ldots,n$ and $i=1,2$) are the vectors in $\R^{2n}$ of the form
\begin{eqnarray*}
e:=\big(e_{11},e_{12},e_{21},e_{22},\ldots,e_{n1},e_{n2}\big)\in
\R^{2n}
\end{eqnarray*}
with 1 at only one position of $e_{ji}$ and $0$ at all the other positions. In order to obtain the optimal solutions to the formulated problem by using Theorem~\ref{Thr7} in what follows, we choose for convenience the norm $\|(x^{j1},x^{j2})\|:=|x^{j1}|+|x^{j2}|$ for each component $x^j\in\R^2$ of $x\in\R^{2n}$.\vspace*{-0.05in}

\section{Optimal Control of Marine Surface Vehicles}\label{marine}

In this section, we explore the dynamic movements of unmanned surface vehicles (USVs) with diverse shapes as they navigate through maritime environments. These USVs are not only capable of reaching their destinations safely but also adept at avoiding both stationary and moving obstacles. This modeling approach transforms the scenario into a controlled perturbed sweeping process using the procedure established in Section~~\ref{model}, which treats each USV as an object within the dynamic system alongside $n-1$ other obstacles on the way to reach the target. This framework finds practical applications in a variety of real-world scenarios, such as oceanographic research, environmental monitoring, maritime surveillance, and autonomous shipping while showcasing the versatility and adaptability of unmanned surface vehicles across different domains.

Note that, in the absence of obstacles, the desired velocities are given by $V(x)=g(x)\in\R^{2n}$. In the presence of obstacles, we take into account the algorithmic developments of \cite{hb1} and seek optimal velocities to escape from different surrounding obstacles by solving the following convex constrained optimization problem:
\begin{eqnarray}\label{P1}
\mbox{minimize }\;&&\|g(x,u)-V(x)\|^2\nonumber\\
\mbox{subject to }\;&&V(x)\in C_{h}(x),
\end{eqnarray}
where the control $u$ is involved into the desired velocity term to adjust the actual velocities of the USVs and make sure that they do not overlap. The algorithmic design in \eqref{P1} means that $V_{k+1}$ is selected as the (unique) element from the set of admissible velocities as the closest one to the desired velocity $g(x,u)$ in order to avoid the overlapping. Consequently, the proposed scheme seeks new directions $V(x)$ of USVs close to the desired direction $g(x,u)$ in order to by pass the surrounding obstacles. The desired position of the next configuration of marine vessels is then generated as $x_{\textrm{ref}}(t+h)=(x_{\textrm{ref}},y_{\textrm{ref}})=x(t)+hV(x)$, and the desired via point posture position of the designated marine craft is given by
$$\eta_{\textrm{ref}}=(x_{\textrm{ref}},y_{\textrm{ref}},\psi_{\textrm{ref}}),\;\mbox{ where }\;
\psi_{\textrm{ref}}:=\tan^{-1}\(\dfrac{y_{\textrm{ref}i}-y_i(t)}{x_{\textrm{ref}i}-x_i(t)}\).
$$

With the constructions in Section~\ref{model}, we now formulate the {\em sweeping optimal control problem} (P) that can be treated as a continuous-time counterpart of the discrete algorithm of the marine surface vehicles model by taking into account the model goals discussed above. Consider the cost functional
\begin{equation}\label{J1}
\mbox{minimize }\;J[x,u,T]:=\disp\frac{1}{2}\big\|x(T)\big\|^2+\dfrac{1}{2}T^2,
\end{equation}
which reflects the model goals to {\em minimize the distance and the time} for the USV  to achieve the target from the admissible configuration set. We describe the continuous-time dynamics by the controlled sweeping process
\begin{equation} \label{st1}
\left\{\begin{array}{lcl}-\dot{x}(t)\in N\big(x(t);C(t)\big)+g\big(x(t),u(t)\big),\\
x(0)=x_0\in C(0),\;u(t)\in U,
\end{array}\right.
\end{equation}
for a.e. $t\in[0,T]$, where the constant set $C(t)$ is taken from \eqref{setC}, the control constraints reduce to \eqref{ut}, and the dynamic nonoverlapping condition $\|x^i(t)-x^j(t)\|\ge R_i+R_j$ is equivalent to the pointwise state constraints
\begin{eqnarray*}
x(t)\in C(t)\Longleftrightarrow\la x^j_*(t),x(t)\ra\le c_j(t)\;\mbox{ for all }\;t\in[0,T]\mbox{ and }\;j=1,\ldots,n-1,
\end{eqnarray*}
which follow from the construction of $C(t)$ and the normal cone definition.

From now on, we exclusively study the sweeping optimal control problem defined in \eqref{J1} and \eqref{st1} with the marine surface vehicles model data specified below. Applying Theorem~\ref{Thr7} allows us to obtain the necessary optimality conditions for the problem that are formulated entirely in terms of the model data. 

Consider the two marine crafts (MCs), MC 1 and MC 2, as moving vehicles. The marine surface vehicles are represented by triangle shapes immersed in discs (Figure~1). The objective is to move MC 1 and MC 2 to the target without colliding with each other. However, in the presence of MC 2, after the contacting time $t^*$ the vehicle MC 1 pushes MC 2 to the target together with the same velocity. The mathematical USV model is taken from the physical ship called Cyber-Ship \cite{kokot}
with the mass of 23.8 kg and the length of 1.255 m. The initial configuration is $x(0)=\(x^1(0),x^2(0)\)$, where $x^1(0)$ is the initial position of MC 1 and $x^2(0)$ is the initial position of MC 2. The position of the target is the origin. The radii of different discs used in this model are chosen as $R_1=R_2=3.5$ m. We specify the other model data as follows:
\begin{equation}\label{sw1}
\left\{
\begin{array}{ll}
n=2,\;x(0)=(x^1(0),x^2(0)),\;x_*(t)=\(1,1,-1,-1\),\;\;c(t)=-7 \mbox{ for all } t,\\[1ex]
g(x,u):=u,\;\vph(x,T):=\dfrac{1}{2}\|x(T)\|^2+\dfrac{1}{2}T^2,\;s_1=s_2=1,\\[1ex]
U=\big\{u(t)=(u^1(t),u^2(t))\in\R^4\;\big|\;\|u^1(t)\|\leq 100,\;\|u^2(t)\|\leq 60\big\},\\[1ex]
\th^u:=\th^u_1=\th^u_2=45^{\circ},\;x^1(0)=\(-25,-25\),\;x^2(0)=\(-15,-15\).
\end{array}\right.
\end{equation}
The moving set $C(t)$ in the sweeping inclusion \eqref{st1} is described now by
$$
\begin{aligned}
C(t)&=\big\{x\in\R^{4}\;\big|\;\langle x_*(t),x\rangle\le c(t)\big\}\nonumber\\
&=\big\{x\in\R^{4}\;\big|\;x^{11}+x^{12}-x^{21}-x^{22}\le -7\big\}\nonumber\\
&=\big\{x\in\R^{4}\;\big|\;|x^{21}-x^{11}|+|x^{22}-x^{12}|\ge 7\big\}\nonumber\\
&\mbox{(under the imposed assumptions}\; x^{21}>x^{11}\mbox{ and }x^{22}>x^{12})\\
&=\big\{x\in\R^{4}\;\big|\;\|x^2-x^1\|\ge 2R\big\} \;\mbox{ for all }\;t\in[0,T].
\end{aligned}
$$
The structure of the problem in \eqref{st1} and \eqref{sw1} suggests that the object changes its velocity only when it hits the boundary at some time $t_m$ with $m\in\{0,1,\ldots,k\}$. If furthermore $t_m<T$, the object slides on the boundary of $C(t)$ for the entire time interval $[t_m,T]$. 

Observing that the controlled sweeping process in \eqref{st1}, \eqref{sw1} satisfies all the assumptions of Theorem~\ref{Thr7}, we have the following relationships to find optimal solutions
$\ox(t) = (\ox^1(t), \ox^2(t))$, $\ou(t) = (\ou^1(t),\ou^2(t))$, and $\oT$.
\begin{figure}[ht]
\centering
\includegraphics[scale=0.36]{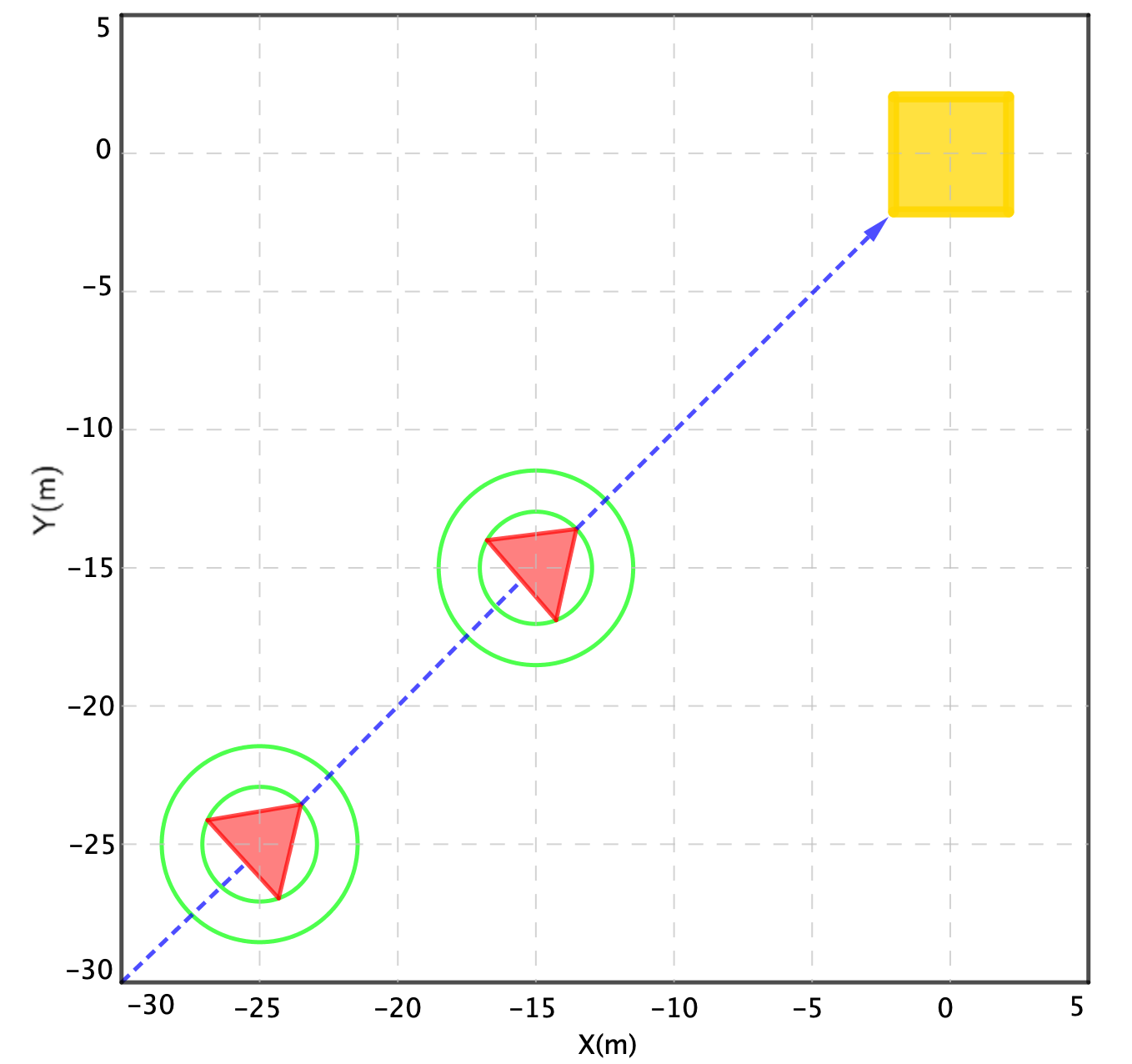}
\caption{Marine vehicle navigation before the contacting time.}
\label{Fig1}
\end{figure}
\begin{figure}[ht]
\centering
\includegraphics[scale=0.36]{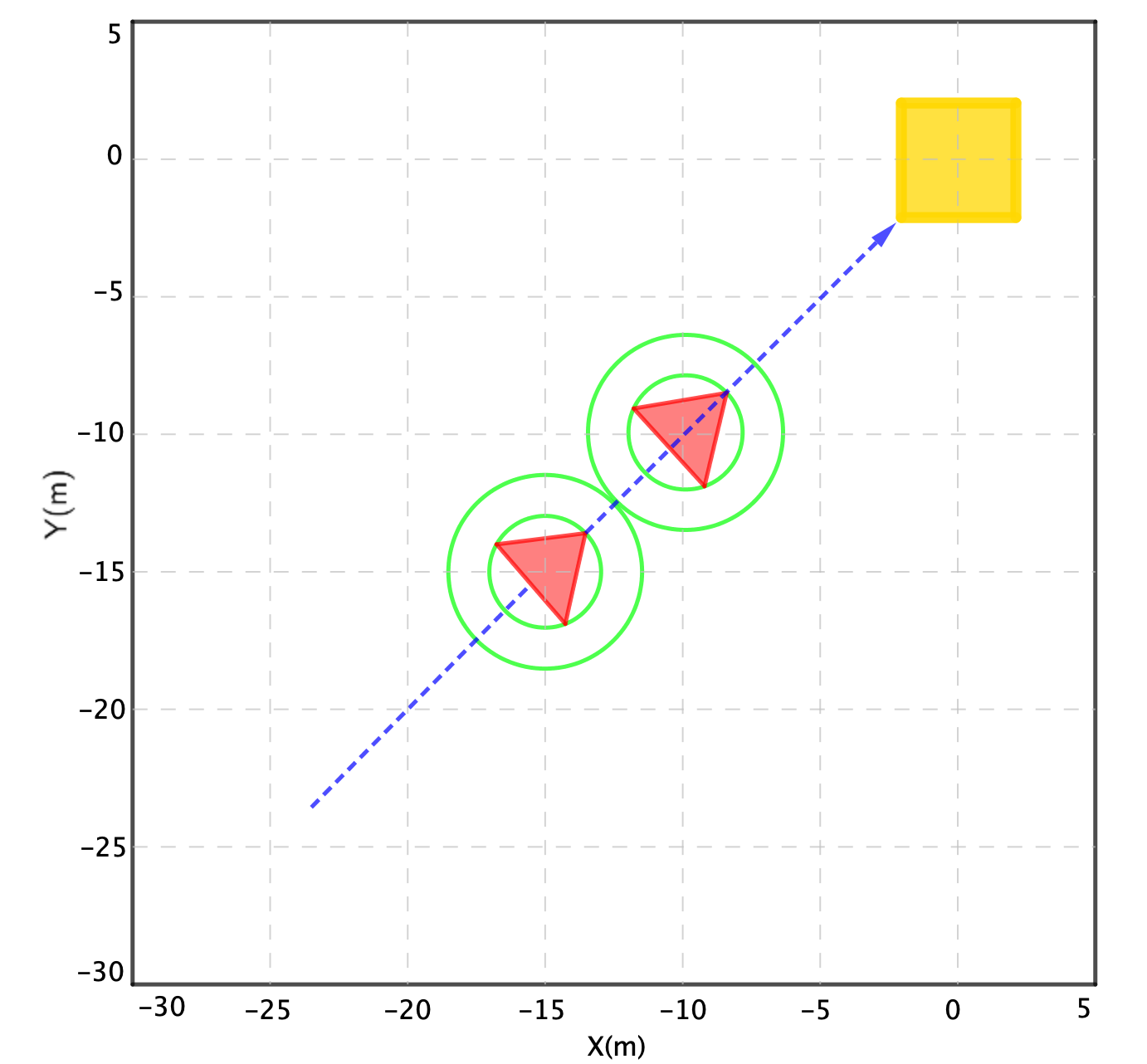}
\caption{Marine vehicle navigation after the contacting time.}
\label{Fig1b}
\end{figure}
\begin{enumerate}
\item $\dot\ox(t)=\nn
\begin{array}{lll} 
\ou(t)-\eta_0(1,1,-1,-1)&\mbox{if}&t\in(0, t_m),\\[2ex]
\ou(t)-\eta_m(1,1,-1,-1)&\mbox{if}&t\in(t_m,\oT)\\[2ex]
\mbox{with }\;\eta_0,\eta_m\ge 0.
\end{array}\right.$
\item $\dot{p}(t)=0$ for a.e. $t\in[0,\oT]$.
\item $q(t)=p(t)-\int^{\oT}_{t}d\gg(\tau)(1,1,-1,-1)$.
\item $\psi(t)=q(t) \in N(\ou(t); U)$ a.e. $t\in [0,\oT]$.
\item $\disp\langle\psi(t),\ou(t)\rangle=\max_{u\in U}\langle\psi(t), u\rangle$ a.e. $t\in [0,\oT]$, which can be rewritten as 
$$\la\psi_1(t),\ou^1(t)\ra+\la\psi_2(t),\ou^2(t)\ra=\max_{(u^1,u^2)\in U}\big\{\la\psi_1(t),u^1\ra+\la\psi_2(t),u^2\ra\big\},$$
where $\psi_1(t)= (\psi_{11}(t),\psi_{12}(t))$ and $\psi_2(t)= (\psi_{21}(t),\psi_{22}(t))$. 
\item $\eta(t)\ge 0,\;\|\ox^2-\ox^1\|>7\Longrightarrow\eta(t)=0$ for a.e. $t\in[0,
\oT]$.
\item $\eta(t)> 0\Longrightarrow\n \ox^2-\ox^1\en=7$ a.e. $t\in[0,\oT]$.
\item
$\nn\begin{array}{lll} 
-p(\oT)=\eta(\oT)(1,1,-1,-1)+\mu\ox(\oT),\;\eta(\oT)\ge 0,\\[2ex]
\eta(\oT)>0 \Longrightarrow \|\ox^1(\oT)+\ox^2(\oT)\|=7,\\[2ex]
\dfrac{1}{\oT}\int^{\oT}_{0}\langle p(t),\dot{\ox}(t)\rangle dt=\mu\oT.
\end{array}
\right.$
\item 
$\left\{\begin{array}{ll} 
(\mu, p, \n\gg_0\en_{TV},\n\gg_>\en_{TV})\neq 0,\;\mu\ge 0,\\
\mu=1\;\mbox{ if }\;\la(1,1,-1,-1),\ox(t)\ra<-7\;\mbox{ on }\;[0,\oT].
\end{array}
\right.$
\item 
$\left\{\begin{array}{ll} 
\mathrm{supp}(\gg_0)\cup\;\mathrm{supp}(\gg_>)\subset
\{t\;|\:\|\ox^1(t)-\ox^2(t)\|=7\},\\
\mathrm{supp}(\gg_>)\cap \mathrm{int}(E_0)=\emptyset\;\textrm{ if }\;\mathrm{int}(E_0)=\mathrm{int}(\{t\;|\;\eta(t)>0\})\ne\emp. 
\end{array}
\right.$
\end{enumerate}
By condition {\bf 2}, $p(t)$ is constant on $[0,\oT]$. Substituting it into {\bf 3} yields
\begin{equation}\label{3'}
q(t)=p(\oT)-\int^{\oT}_{t}(1,1,-1,-1)d\gg(\tau)=p(\oT)-(1,1,-1,-1)\gg([t,\oT]).
\end{equation}
Then combining \eqref{3'} and {\bf 8} gives us the expression
$$
q(t)=-\mu\ox(\oT)-(1,1,-1,-1)\big(\eta(\oT)+\gg([t,\oT])\big)\;\mbox{ on }\;[0,\oT].
$$
Since we start from an interior point of the polyhedron, the measure $\gg$ is zero in the interval $[0,t_m]$. Let us investigate behavior of the trajectory before and after the hitting time $t=t_m$. Remembering that the trajectory does not hit the boundary before $t=t_m$, it follows from  {\bf 6} that $\gg([t,T])\equiv\gg([t_m,\oT])$ for all $t\in[0,t_m)$. Thus 
\begin{equation}\label{r_c}
q(t)\equiv-\mu\ox(\oT)-(1,1,-1,-1)\big(\eta(\oT)+\gg([t_m,\oT])\big)=:r_m=(r_{1m},r_{2m})\in\R^4
\end{equation}
for all $t\in[0,t_m)$. From now on, we suppose for determinacy and simplicity in this complicated model that optimal control under consideration is constant on $[0,\oT]$ and belong to the interior of $U$. Combining {\bf 4} and {\bf 5} tells us that
$$
\la r_{1m},\ou^1\ra+\la r_{2m},\ou^2\ra=\max_{(u^1,u^2)\in U}\big\{\la r_{1m},u^1\ra+\la r_{2m},u^2\ra\big\},
$$
which can be written by \eqref{r_c} in the form
$$
\begin{aligned}
&\big\la-\mu \ox^1(\oT)-(1,1)(\eta(\oT)+\gg([t_m,\oT])),\ou^1\big\ra +
\big\la-\mu \ox^2(\oT)-(-1,-1)(\eta(\oT)+\gg([t_m,\oT])),\ou^2\big\ra\\
&= \max_{(u^1,u^2)\in U}\big\{\big\la-\mu \ox^1(\oT)-(1,1)(\eta(\oT)+\gg([t_m,\oT])),u^1\big\ra +\big\la-\mu \ox^2(\oT)-(-1,-1)(\eta(\oT)+\gg([t_m,\oT])),u^2\big\ra\big\}
\end{aligned}
$$
According to the above, consider the linear function
$$\phi(u^1,u^2):= \big\la-\mu \ox^1(\oT)-(1,1)(\eta(\oT)+\gg([t_m,\oT)),u^1\big\ra +\big\la-\mu \ox^2(\oT)-(-1,-1)(\eta(\oT)+\gg([t_m,\oT])),u^2\big\ra, $$
which attains its maximum over $U$ at the interior point $(\ou_1,\ou_2)$. Thus the Fermat rule for $\phi$ gives us the conditions
$$-\mu \ox^1(\oT)=(1,1)(\eta(\oT)+\gg([t_m,\oT]))\;\mbox{ and }\;
\mu \ox^2(\oT)=(1,1)(\eta(\oT)+\gg([t_m,\oT])),
$$
which imply in turn the equalities
\begin{equation}\label{u-T}
-\ox^1(\oT)= \ox^2(\oT)=(1,1)(\eta(\oT)+\gg([t_m,\oT])).
\end{equation}
When the particle moving hits the boundary at the time $ t=t_m$, it follows from {\bf 6} that 
$$|\ox^{21}(t_m)-\ox^{11}(t_m)|+|\ox^{22}(t_m)-\ox^{12}(t_m)|=7,\quad\mbox{and}
$$
$$
\eta(t)=\begin{cases}
\eta_0 = 0 \;\;\mbox{if}\;\;t\in [0,t_m],\\
\eta_m\ge 0\;\;\mbox{if}\;\;t\in[t_m,\oT].
\end{cases}
$$
After the hitting point, both these marine USVs will move to the target under the control $\ou$. This yields
\begin{eqnarray*}
\left\{\begin{array}{ll}
\ox^1(t)=\(-25+t\|\ou^{1}\|\cos\theta^u,-25+t\|\ou^{1}\|\sin\theta^u\),\\
\ox^2(t)=\(-15+t\|\ou^{2}\|\cos\theta^u,-15+t\|\ou^{2}\|\sin\theta^u\)
\end{array}\right.
\end{eqnarray*}
for $t\in[0,t_m)$ as well as the representations
\begin{eqnarray*}
\left\{\begin{array}{ll}
\ox^1(t)=\(-25+t_m\|\ou^{1}\|\cos\theta^u+(t-t_m)(\|\ou^{1}\|\cos\th^u-\eta_m),-25+t_m\|\ou^{1}\|\sin\theta^u+(t-t_m)(\|\ou^{1}\|\sin\th^u-\eta_m)\),\\
\ox^2(t)=\(-15+t_m\|\ou^{2}\|\cos\theta^u+(t-t_m)(\|\ou^{2}\|\cos\th^u+\eta_m),-15+t_m\|\ou^{2}\|\sin\theta^u+(t-t_m)(\|\ou^{2}\|\sin\th^u+\eta_m)\)
\end{array}\right.
\end{eqnarray*}
for $t\in[t_m,\oT]$.
Recalling that the normal vectors are inactive before the hitting time, we arrive at the equality
$$
2\Big|10 +\dfrac{\sqrt{2}}{2} t_m(\|\ou^{2}\|-\|\ou^{1}\|)+\(t-t_m\)\Big(\dfrac{\sqrt{2}}{2}( \|\ou^{2}\|-\|\ou^{1}\|)+2\eta_m\Big)\Big|=7\;\mbox{ if }\;t\ge t_m
$$
for the time of contact. This brings us to the equation
$$\Big|10 -\dfrac{\sqrt{2}}{2}t_m(\|\ou^{2}\|-\|\ou^{1}\|)\Big|=7,$$
which has the following two solutions:
\begin{equation}\label{t_m}
t_m=\dfrac{-13}{\sqrt{2}(\|\ou^{2}\|-\|\ou^{1}\|)}\;\mbox{ and }\;t_m= \dfrac{-27}{\sqrt{2}(\|\ou^{2}\|-\|\ou^{1}\|)}
\end{equation}
with $\|\ou^2\|<\|\ou^1\|$, which will both be checked for the optimality below. When $\ox(\cdot)$ hits the boundary of $C$, it would stay there while pointing in the direction shown in Figure~2. Thus we get
\begin{equation}\label{u-eta}
\dfrac{\sqrt{2}}{2}(\|\ou^{2}\|-\|\ou^{1}\|)+2\eta_m=0,\;\mbox{ and so }\;\eta_m=\dfrac{\sqrt{2}(\|\ou^{1}\|-\|\ou^{2}\|)}{4}.
\end{equation}
Furthermore, it follows from \eqref{u-T} that the optimal time is calculated by
$$\oT=\dfrac{40\sqrt{2}}{\|\ou^{2}\|+\|\ou^{1}\|}.$$
Since $p(\cdot)$ is constant on $[0,\oT]$, the third condition in {\bf 8} yields
\begin{equation}\label{p}
\dfrac{1}{\oT}\int^{\oT}_{0}\langle p(\oT),\dot{\ox}(t)\rangle dt=\mu\oT,\; \;\mbox{i.e.,}\;\;p(\oT)\ox(\oT)=\mu{\oT}^2.
\end{equation}
To get enough information, suppose that $\mu=1$. Substituting \eqref{p} into the first equation in {\bf 8} yields 
\begin{equation}\label{T}
-\eta_m(1,1,-1,-1)=\ox(\oT)+\frac{\oT^2}{\ox(\oT)}.
\end{equation}
With \eqref{t_m} and $\eta_m$ taken from \eqref{t_m} and \eqref{u-eta}, respectively, the cost function value at \eqref{T} is 
$$
\begin{array}{ll}
J[\ox,\ou,\oT]&=\[-25+\dfrac{\sqrt{2}}{2} t_m\|\ou^{1}\|+(t-t_m)\(\dfrac{\sqrt{2}}{2}\|\ou^{1}\|-\eta_m\)\]^2\\
&+\[-15+\dfrac{\sqrt{2}}{2} t_m\|\ou^{2}\|+(t-t_m)\(\dfrac{\sqrt{2}}{2}\|\ou^{2}\|+\eta_m\)\]^2 +\disp\frac{1}{2}\oT^2.
\end{array}
$$

When the object hits the boundary, it then slights there until the end of the process, i.e., $t\ge t_m$. It then follows that $\eta_m>0$.
By using computer calculations, all the cases can be analyzed.
This shows that the best performance is obtained by choosing the contacting time $t_m=\dfrac{-13}{\sqrt{2}(\|\ou^{2}\|-\|\ou^{1}\|)}$ for which we have
$$
\begin{aligned}
J[\ox,\ou,\oT]=& \Big[-25+\dfrac{-13\|\ou^1|}{2(\|\ou^{2}\|-\|\ou^{1}\|)}+\Big(\dfrac{40\sqrt{2}}{\|\ou^{2}\|+\|\ou^{1}\|}-\dfrac{-13}{\sqrt{2}(\|\ou^{2}\|-\|\ou^{1}\|)}\Big)\Big(\dfrac{\sqrt{2}}{2}\|\ou^{1}\|-\dfrac{\sqrt{2}(\|\ou^{1}\|-\|\ou^{2}\|)}{4}\Big)\Big]^2\\
&+\Big[-15+\dfrac{-13\|\ou^2\|}{2(\|\ou^{2}\|-\|\ou^{1}\|)}+\Big(\dfrac{40\sqrt{2}}{\|\ou^{2}\|+\|\ou^{1}\|}-\dfrac{-13}{\sqrt{2}(\|\ou^{2}\|-\|\ou^{1}\|)}\Big)\Big(\dfrac{\sqrt{2}}{2}\|\ou^{2}\|+\dfrac{\sqrt{2}(\|\ou^{1}\|-\|\ou^{2}\|)}{4}\Big)\Big]^2\\
&+\dfrac{1}{2}\Big(\dfrac{40\sqrt{2}}{\|\ou^{2}\|+\|\ou^{1}\|}\Big)^2.
\end{aligned}
$$ 
This function achieves the minimum value $J\approx 6.1875$ at $\|\ou^{1}\|\approx 99.9999$, $\|\ou^{2}\|\approx 59.9999$, and $\oT\approx 0.35355$, which imply that $t_m\approx 0.2298$ and $\ox(\oT)\approx (-1.7502,-1.7502,1.7498,1.7498)$. This tells us that under the imposed assumptions we arrive at the unique optimal control
$(\ou^1,\ou^2)\approx(70.7106, 70.7106, 42.4263, 42.4263)$. The corresponding optimal trajectory is
\begin{eqnarray*}
\left\{\begin{array}{ll}
\ox^1(t)\approx\(-25+70.7106t,-25+70.7106t\),\\
\ox^2(t)\approx\(-15+42.4263t,-15+42.4263t\)
\end{array}\right.
\end{eqnarray*}
for $t\in[0,0.2298)$, and finally
\begin{eqnarray*}
\left\{\begin{array}{ll}
\ox^1(t)\approx\(-21.7500+56.5685t,-21.7500+56.5685t\),\\
\ox^2(t)\approx\( -18.2500+56.5685t, -18.2500+56.5685t\)
\end{array}\right.
\end{eqnarray*}
for the last interval $t\in[0.2298,0.35355]$.\vspace*{-0.05in}

\section{Applications to Nanoparticles}\label{sec:example2}
\setcounter{equation}{0}\vspace*{-0.1in}

The motivation to explore the motion of nanoparticles in a straight tube is driven by the profound influence of nanotechnology and biotechnology on pharmacology. This impact has significantly enhanced the performance of existing drugs while also enabling the development and utilization of novel drugs and therapies. Recognizing the transformative effects of these technologies on drug delivery and treatment strategies, this section seeks to understand and optimize the motion of nanoparticles within a confined space. Our research is motivated by the potential applications of such motion in improving drug delivery precision, facilitating controlled release mechanisms, and ultimately contributing to advancements in pharmaceutical and therapeutic practices. We study a scenario involving nanoparticles characterized as inelastic disks with varying radii. The objective for each nanoparticle is to traverse straight circular capillaries and reach the target at the end of the tube, in the shortest time. Throughout this motion, the nanoparticles potentially make contact with the other $n-1$ nanoparticles and treat them as obstacles while avoiding collisions. 

For simplicity, taking into account the assumptions in \cite{mb} that the motion of the nanoparticles is considered on the $xz$-plane, i.e., the nanoparticle inelastic disk center is always in the $y = 0$ plane, and there is no rotation of the nanoparticle about the $z$-axis. The nanoparticle therefore has two degrees of freedom: its motion is fully described by specifying $(x^{i1},x^{i2})$ as functions of time $t$, where $x^{i1}$ and $x^{i2}$ are coordinates of the center of nanoparticle $i$. The trajectory of each nanoparticle is governed by the forces exerted by the flow (blood stream) and the gravitation. Forces acting on nanoparticles include hemodynamic forces, buoyancy, Van der Waals' interactions between nanoparticles each-to-other and between nanoparticles and walls. To derive the equations describing the motion on the $x$-axis and on the $z$-axis, we apply (as suggested in \cite{Zhao2}) the balance principle of forces acting on the nanoparticles. As considered in \cite{Zhao1}, the balance of forces acting on the nanoparticles requires, for each nanoparticle $i$, a system of two equations on the $x$-axis and the $z$-axis; see Figure~3.

Since all the nanoparticles intend to reach the end of the straight circular capillaries by the {\em shortest path}, their desired spontaneous (i.e., in the absence of other obstacles) velocities can be represented as
\begin{equation*}
S(x)=\big(s_1,\ldots,s_n\big),\;\mbox{ with }\; s_i=(s_{x^{i1}},s_{x^{i2}}),\;\mbox{ for }i =1,\ldots, n,
\end{equation*}
where the velocity of nanoparticles are obtained after solving the following ODE system from \cite{mb}:
\begin{eqnarray*}
\left\{\begin{matrix}
\dfrac{dx^{i1}}{dt}&=&\dfrac{-\left [ h_6^{w_1}(x^{i2})+ h_6^{w_2}(x^{i2}) \right ]S}{h_1^{w_1}(x^{i2})+ h_1^{w_2}(x^{i2}) }\\ 
\dfrac{dx^{i2}}{dt}&=&\dfrac{\frac{4\pi}{3}gR\Delta \rho +\sum_{j\neq i}F_{ij}+F_{iw}}{h_5^{w_1}(x^{i2})+ h_5^{w_2}(x^{i2})}
\end{matrix}\right. \;\mbox{ for }i =1,\ldots, n.
\end{eqnarray*}
Here $S$ is the shear rate of the flow, $h^{w_k}_1$ and $h^{w_k}_6$ are the hemodynamic resistant force induced by the wall $w_k$ on the $x$-axis, $h^{w_k}_5$ is the hemodynamic resistant force induced by the wall $w_k$ on the $z$-axis ($k= 1,2$), $F_{ij}$ is the force exerted over nanoparticle $i$ due to van der Waals' interaction energy between the nanoparticle $i$ and nanoparticle $j$, and $F_{iw}$ is the force exerted over nanoparticle $i$ due to van der Waals' interaction energy between the nanoparticle $i$ and the walls, i.e.,
$$F_{iw}=|F_{iw_1}-F_{iw_2}|,$$
where $F_{iw_k}$ is the force exerted over nanoparticle $i$ due to van der Waals' interaction energy between the nanoparticle $i$ and the wall $k$. It is known that 
$$F_{ij}=\dfrac{AR^3_i}{D^4_{ij}},$$
where $D_{ij}$ is the distance between the nanoparticles $i$ and $j$. In the equations above, $R_i$ denotes the radius of nanoparticle $i$, $\Delta_\rho$ stands for the difference between mass densities of the nanoparticles and the fluid, and $g$ is the standard gravitational acceleration.

Now we introduce the sweeping optimal control problem $(P)$ for the nanoparticles model discussed above. This formulation considers the aforementioned model objective and utilizes the modeling framework developed in Section~\ref{model}. Similarly to Section~\ref{marine}, define the cost functional
\begin{equation}\label{t:102*}
\mbox{minimize }\;J[x,u,T]:=\disp\frac{1}{2}\big\|x(T)\big\|^2+\disp\frac{1}{2}T^2,
\end{equation}
which reflects the model goals to {\em minimize the distance} of the nanoparticle from the admissible configuration set to the end of the straight circular capillaries together with {\em minimizing the time} to reach the target. Employing the nonoverlapping algorithm introduced in \cite{maury} for the
study of inelastic collision (see also \cite{bk,venel} for further developments in other models) and taking into account the above discussions, we describe the continuous-time dynamics by the {\em controlled sweeping process}
\begin{equation}\label{t:101*}
\left\{\begin{matrix}
\dot{x}(t)\in-N\big(x(t);C(t)\big)+g\big(x(t),u(t)\big)\;\textrm{ a.e. }\;t\in[0,T],\;x(0)=x_0\in C(0)\subset\R^{2n},\\
u(t)\in U\;\textrm{ a.e. }\;t\in[0,T],
\end{matrix}\right.
\end{equation}
where the control constraints reduce to \eqref{ut}, and the dynamic nonoverlapping condition $\|x^i(t)-x^j(t)\|\ge R_i+R_j$ is used to formulate the moving set 
\begin{equation}\label{e:131***}
C(t):=\big\{x(t)\in\R^{2n}\;\big|\;\la x^j_*(t),x(t)\ra\le c_j(t),\;j=1,\ldots,n-1\big\}
\end{equation}
with $c_j(t):=-(R_{nano}+R_{obs})$, where $R_{nano}$, $R_{obs}$ are the radii of the considered nanoparticle and the obstacle, respectively. The $n-1$ vertices of the polyhedron are constructed as in \eqref{e}.

Observe that, although the frameworks of nanoparticle and marine  USV models are rather similar, these two models are essentially different. The difference between them lie in both dynamics and constraints. In particular, the motion of nanoparticles is determined by the forces applied by the flow (bloodstream) and gravity. The forces influencing nanoparticles encompass hemodynamic forces, buoyancy, Van der Waals' interactions among nanoparticles, both inter-particle and with surrounding walls. Consequently, to have the velocities of the nanoparticles, we should take into account all these factors, which are presented in Table~1.

Our initial data for the nanoparticle model are given as follows:
\begin{equation}\label{datanano}
\left\{
\begin{array}{ll}
n=2,\;\;l = 7 \mu m=7000 nm,\;\;h = 0.0075,\;\;R_1=5 nm,\;\;R_2=10 nm,\;\;S=1600e6\; s^{-1},\\[1ex]
A = 5e-21,\;\;g = 9.81e9,\;\;\mu = 3.4e-3,\;\;\Delta\rho=-1e3,\;\;a = 0.5 \mu m=500 nm,\\[1ex]
x_*=\(1,1,-1,-1\),\;\;x(0)=(x^{11}(0),x^{12}(0),x^{21}(0),x^{22}(0)),\;\;c=-R_1-R_2=-15,\\[1ex]
g(x,u):=u,\;\;\vph(x,T):=\dfrac{1}{2}\|x(T)\|^2+\disp\frac{1}{2}T^2,\\[1ex]
U:=\big\{u=(u^1,u^2)\in \R^4|\; \|u^1\|\in [0,3], \|u^2\|\in [0,3]\big|\;\|u^1\|=2\|u^2\|\big\},\\[1ex]
\th^u=45^{\circ},\;x^1(0)=\(-350,-350\),\;x^2(0)=\(-200,-200\).
\end{array}\right.
\end{equation}

\begin{figure}[h!]
\centering
\includegraphics[scale=0.45]{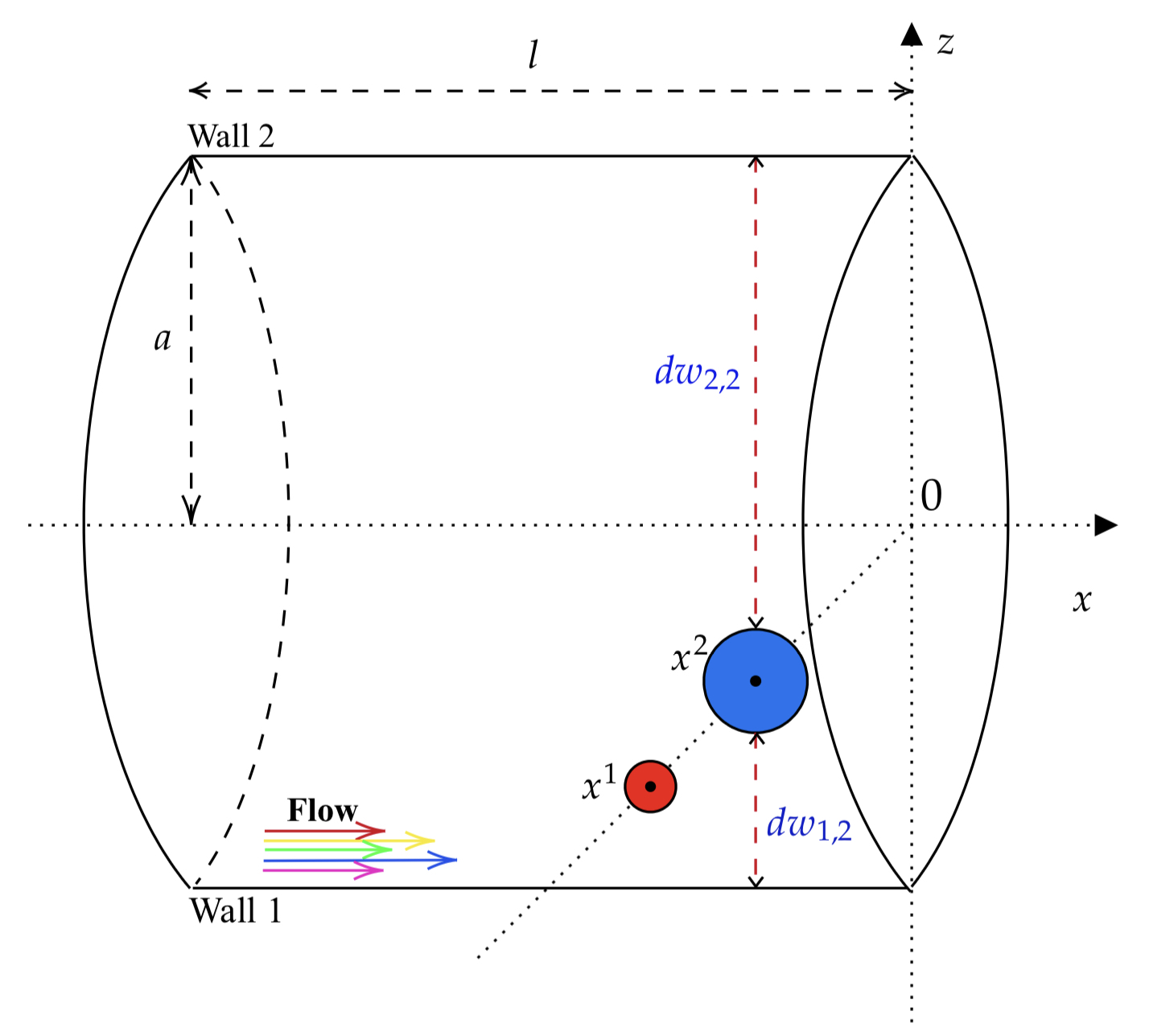}
\caption{Nanoparticles in a straight tube.}
\label{Fig3}
\end{figure}
\begin{figure}[ht]
\centering
\includegraphics[scale=0.45]{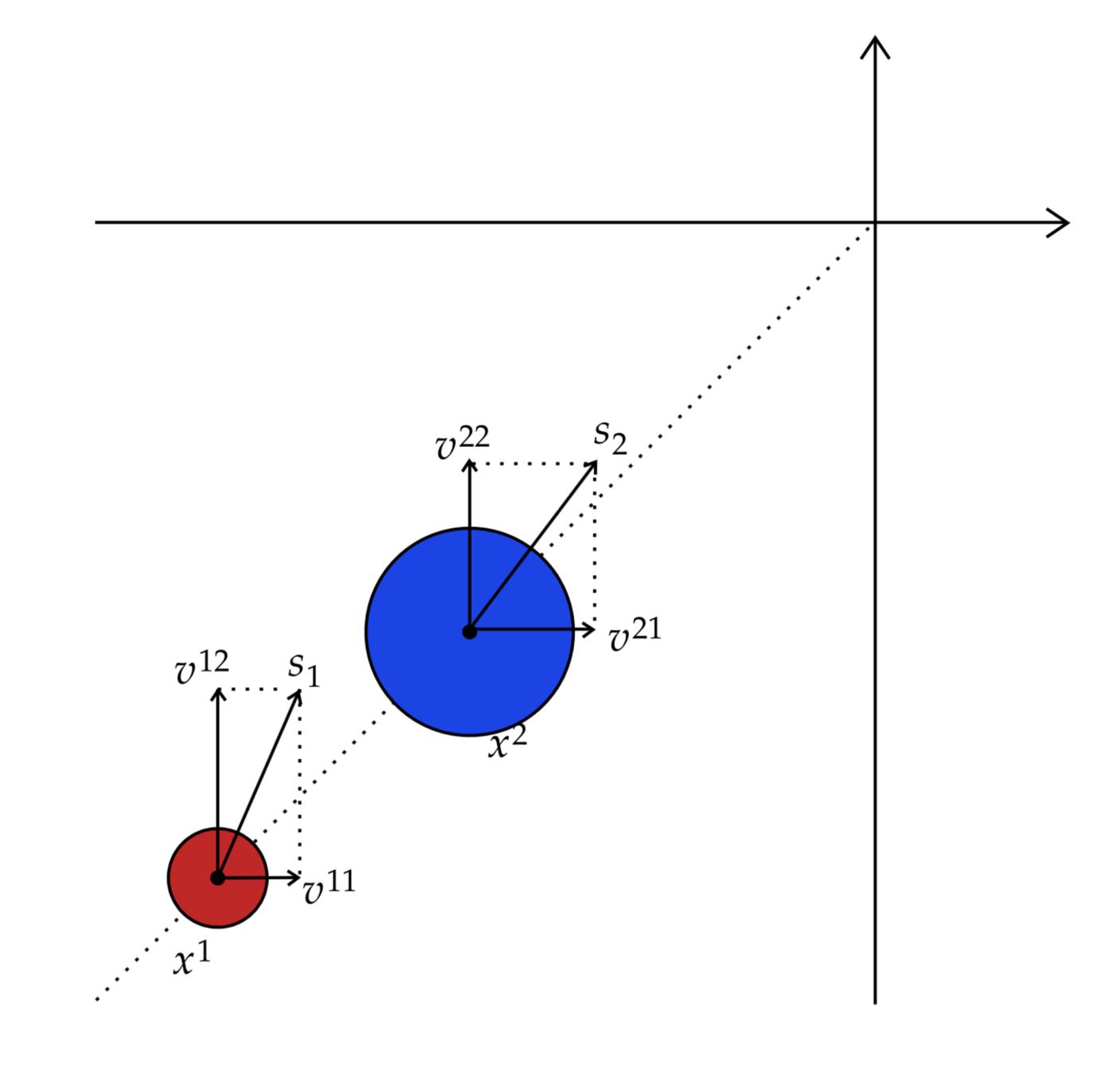}
\caption{Nanoparticles with the desired velocities.}
\label{Fig1a}
\end{figure}

The data for nanoparticles in the viscous flow inside a tube with length $l$, width $a$, and parameters $h$, $R_1$, $R_2$, $T$, $A$, $g$, $\mu$, and $\Delta_p$ are taken from \cite{mb}. The data  for $S$ and related parameters reflected in Figure~4 are taken from \cite{Saka}. The controlled desired velocities are naturally described by 
\begin{eqnarray*}
g\big(x(t),u(t)\big)=u(t)=\big(s_1\|u^1(t)\|\cos\th^u(t),s_1\|u^1(t)\|\sin\th^u(t),s_2\|u^2(t)\|\cos\th^u(t),s_2\|u^2(t)\|\sin\th^u(t)\big).
\end{eqnarray*}
Table~1 collects the results of calculations of the model ingredients according to the above formulas with the data of \eqref{datanano} for two nanoparticles $i=1,2$.

\begin{table}[h]\label{tab1}
\caption{Data of the nanoparticles inside a straight microtube} 
\centering 
\begin{tabular}{|c | c | c |} 
\hline\hline
Nanoparticle $i$ & Data \\ [0.5ex] 
\hline 
$i=1$ & $dw_{1,1}=a-|x^{12}|-R_1=145$  \\  [1ex]
& $dw_{2,1}=a+|x^{12}|-R_1=845$ \\  [1ex]
& $h_6^{w_1}(x^{12})=6\pi\mu R_1\left [ 1+\frac{9}{16}\left (\frac{R_1}{R_1+dw_{1,1}} \right )+0.13868\left ( \frac{R_1}{R_1+dw_{1,1}} \right )^{1.4829} \right ]\approx 0.3267374$ \\ [1ex]
& $h_6^{w_2}(x^{12})=6\pi\mu R_1\left [ 1+\frac{9}{16}\left ( \frac{R_1}{R_1+dw_{2,1}} \right )+0.13868\left ( \frac{R_1}{R_1+dw_{2,1}} \right )^{1.4829} \right ]\approx 0.3215246$ \\  [1ex]
 & $h_1^{w_1}(x^{12})=\frac{48}{15}\pi\mu R_1\left [ ln\left ( 1-\frac{R_1}{R_1+dw_{1,1}} \right ) +\frac{R_1}{R_1+dw_{1,1}}\right ]\approx -0.0000971$ \\ [1ex] 
 &$h_1^{w_2}(x^{12})=\frac{48}{15}\pi\mu R_1\left [ ln\left ( 1-\frac{R_1}{R_1+dw_{2,1}} \right ) +\frac{R_1}{R_1+dw_{2,1}}\right ]\approx -0.0000030$   \\ [1ex]
 & $h_5^{w_1}(x^{12})=h_5^{w_2}(x^{12})=6\pi\mu R_1 \left [ -\frac{1}{1-\frac{R_1}{|x^{12}|}}+\frac{1}{8}ln\left ( 1-\frac{R_1}{|x^{12}|} \right ) \right ]\approx -0.3256629 $ \\ [1ex]
 & $F_{1w}= \Big|\dfrac{AR^3_1}{dw_{1,1}^4}-\dfrac{AR^3_1}{dw_{2,1}^4}\Big|\approx 1.325e-27$ \\ [1ex]
 & $F_{12}=\dfrac{AR^3_1}{D^4_{12}}\approx 4.139e-28,\textrm{ where }D_{12}= 197.1320344
$ \\ [1ex]
 & $
\left\{\begin{matrix}
v^{11}=\dfrac{dx^{11}}{dt}&=&\dfrac{-\left [ h_6^{w_1}(x^{12})+ h_6^{w_2}(x^{12}) \right ]S}{h_1^{w_1}(x^{12})+ h_1^{w_2}(x^{12}) }\approx 1.036e+13\\ 
v^{12}=\dfrac{dx^{12}}{dt}&=&\dfrac{\frac{4\pi}{3}gR_1\Delta \rho +F_{12}+F_{1w}}{h_5^{w_1}(x^{12})+ h_5^{w_2}(x^{12})}\approx 3.154e+14
\end{matrix}\right.,
$ \\ [3ex]
&$s_1=3.156e+14$ \\ [2ex]
$i=2$ & $dw_{1,2}=a-|x^{22}|-R_2=290$  \\  [1ex]
& $dw_{2,2}=a+|x^{22}|-R_2=690$ \\  [1ex]
& $h_6^{w_1}(x^{22})=6\pi\mu R_2\left [ 1+\frac{9}{16}\left (\frac{R_2}{R_2+dw_{1,2}} \right )+0.13868\left ( \frac{R_2}{R_2+dw_{1,2}} \right )^{1.4829} \right ]\approx  0.6534748$\\ [1ex]
& $h_6^{w_2}(x^{22})=6\pi\mu R_2\left [ 1+\frac{9}{16}\left ( \frac{R_2}{R_2+dw_{2,2}} \right )+0.13868\left ( \frac{R_2}{R_2+dw_{2,2}} \right )^{1.4829} \right ]\approx 0.6461981$ \\  [1ex]
& $h_1^{w_1}(x^{22})=\frac{48}{15}\pi\mu R_2\left [ ln\left ( 1-\frac{R_2}{R_2+dw_{1,2}} \right ) +\frac{R_2}{R_2+dw_{1,2}}\right ]\approx -0.0001942$\\  [1ex]
& $h_1^{w_2}(x^{22})=\frac{48}{15}\pi\mu R_2\left [ ln\left ( 1-\frac{R_2}{R_2+dw_{2,2}} \right ) +\frac{R_2}{R_2+dw_{2,2}}\right ]\approx -0.0000352$ \\  [1ex]
& $h_5^{w_1}(x^{22})=h_5^{w_2}(x^{22})=6\pi\mu R_2 \left [ \frac{1}{1-\frac{R_2}{|x^{22}|}}+\frac{1}{8}ln\left ( 1-\frac{R_2}{|x^{22}|} \right ) \right ]\approx -0.6787248 $  \\ [1ex]
& $F_{2w}= \Big|\dfrac{AR^3_2}{dw_{1,2}^4}-\dfrac{AR^3_2}{dw_{2,2}^4}\Big|\approx 6.849e-28$ \\ [1ex]
& $F_{21}=\dfrac{AR^3_2}{D^4_{21}}\approx 3.311e-27,\textrm{ where }D_{21}= 197.1320344$\\ [1ex]
& $
\left\{\begin{matrix}
v^{21}=\dfrac{dx^{21}}{dt}&=&\dfrac{-\left [ h_6^{w_1}(x^{22})+ h_6^{w_2}(x^{22}) \right ]S}{h_1^{w_1}(x^{22})+ h_1^{w_2}(x^{22}) }\approx 9.064e+12\\ 
v^{22}=\dfrac{dx^{22}}{dt}&=&\dfrac{\frac{4\pi}{3}gR_2\Delta \rho +F_{21}+F_{2w}}{h_5^{w_1}(x^{22})+ h_5^{w_2}(x^{22})}\approx 3.027e+14
\end{matrix}\right.,
$ \\ [2ex]
 & $s_2=3.0284e+14$ \\ [2ex]
\hline 
\end{tabular}
\end{table}

At the initial time, we have $(x^{11}(0)-x^{21}(0))^2+(x^{12}(0)-x^{22}(0))^2=45000$. Let $t_\ast>0$ stand for the first contacting time between the two nanoparticles, i.e.,
\begin{eqnarray*}
t_\ast:=\min\big\{t\in[0,T]\;\big|\;\|\ox^1(t)-\ox^2(t)\|=R_1+R_2\big\}.
\end{eqnarray*}
Recall that the nanoparticle tends to keep its constant direction and velocity until either touching the other nanoparticle (obstacle), or reaching the end of the straight microtube at $t=\oT$. To proceed further, suppose that $\mu > 0$ for definiteness and apply the necessary optimality conditions of Theorem~\ref{Thr7} in our control model. We arrive at the following relationships:
\begin{itemize}
\item[\bf(1)] $-\big(\dot{\ox}^{11}(t),\;\dot{\ox}^{12}(t),\;\dot{\ox}^{21}(t),\;\dot{\ox}^{22}(t)\big)=\eta(t)(1,1,-1,-1)-\\
\big(s_1|\ou^1(t)|\cos\th^u,\;s_1|\ou^1(t)|\sin\th^u,\;s_2|\ou^2(t)|\cos\th^u,\;s_2|\ou^2(t)|\sin\th^u\big)$ for a.e.\ $t\in[0,\oT]$.
\item[\bf(2)] $\|\ox^{1}(t)-\ox^{2}(t)\|>R_1+R_2\Longrightarrow\eta(t)=0$ for a.e.\ $t\in[0,\oT]$.
\item[\bf(3)] $\eta(t)>0\Longrightarrow\la x_*,q(t)\ra=0$ for a.e.\ $t\in[0,\oT]$.
\item[\bf(4)] $
(-p(\oT)-\eta(\oT) x_\ast,\bar{H})\in\nabla\varphi(\bar{x}(\oT),\oT)\;
$
where $\bar{H}:=\dfrac{1}{\oT}\int_{0}^{\oT}\langle p(t),\dot{\bar{x}}(t)\rangle dt$.
\item[\bf(5)] $q(t)=p(t)-\disp\int_{(t,\oT]}d\gg\(\tau\)x_*$ for all $t\in[0,\oT]$ except at most a countable subset.
\item[\bf(6)] $\big\la \psi(t),\ou(t)\big\ra=\max_{u\in U}\big\la \psi(t),u\big\ra,\;\mbox{ where }\;\psi(t):=\nabla_u g\big(\ox(t),\ou(t)\big)^*q(t)\mbox{ for a.e. }\;t\in[0,\oT]$.
\item[\bf(7)] $\eta(\oT)>0 \Longrightarrow \|\ox^1(\oT)+\ox^2(\oT)\|=R_1+R_2.$
\item[\bf(8)]$\dot{p}(t)=0,\;\mbox{ for a.e. }\;t\in[0,\oT]$
\end{itemize}
Taking into account that the nanoparticle directions are constant as well as the assumptions in the model imposed above, we seek for simplicity constant optimal controls. Then it follows from condition {\bf 2} that the function $\eta(\cdot)$ is piecewise constant on $[0,\oT]$ and admits the representation
\begin{eqnarray}\label{eta}
\eta(t)=\left\{\begin{array}{ll}
0\quad\mbox{ for }\;t\in[0,t_\ast),\\
\eta\quad\mbox{ for }\;t\in[t_\ast,\oT].
\end{array}\right.
\end{eqnarray}
Using now {\bf 1}, the dynamic equations prior to and after the time $t_\ast$ can be rewritten as
\begin{eqnarray*}
\left\{\begin{array}{ll}
\dot{\ox}^1(t)=\(s_1|\ou^1|\cos\th^u,s_1|\ou^1|\sin\th^u\)\;\mbox{ and}\\
\dot{\ox}^2(t)=\(s_2|\ou^2|\cos\th^u,s_2|\ou^2|\sin\th^u\)\;\mbox{ for }\;t\in[0,t_\ast),
\end{array}\right.
\end{eqnarray*}\vspace*{-0.2in}
\begin{eqnarray*}
\left\{\begin{array}{ll}
\dot{\ox}^1(t)=\(-\eta(t)+s_1|\ou^1|\cos\th^u,-\eta(t)+s_1|\ou^1|\sin\th^u\)\;\mbox{ and}\\
\dot{\ox}^2(t)=\(\eta(t)+s_2|\ou^2|\cos\th^u,\eta(t)+s_2|\ou^2|\sin\th^u\)\;\mbox{ for }\;t\in[t_\ast,\oT].
\end{array}\right.
\end{eqnarray*}
At the contacting time $t=t_{\ast}$, the two nanoparticles have the same velocities till the end of the straight tube, which yields $\dot{\ox}^1(t)=\dot{\ox}^2(t)$ for all $t\in[t_\ast,\oT]$. This gives us in turn the following calculation of the corresponding value of $\eta$ in representation \eqref{eta}: 
\begin{eqnarray}\label{eq-eta2}
\eta=\left\{\begin{array}{ll}
\frac{1}{2}\(s_1|\ou^1|\cos\th^u-s_2|\ou^2|\cos\th^u\)\;&\mbox{ if }\;s_1|\ou^1|\ne  s_2|\ou^2|\;\mbox{ and }\;\cos\th^u=\sin\th^u,\\
0&\mbox{ otherwise}
\end{array}\right.
\end{eqnarray}
The case of $\eta=0$ is trivial. Applying the first formula in \eqref{eq-eta2} and the constant nanoparticles velocity after the contacting time till reaching the target, we deduce from the Newton-Leibniz formula that
\begin{eqnarray*}
\left\{\begin{array}{ll}
\ox^1(t)=\(\ox^{11}(0),\ox^{12}(0)\)+\(ts_1|\ou^1|\cos\th^u,ts_1|\ou^1|\sin\th^u\)\;\mbox{ and}\\
\ox^2(t)=\(\ox^{21}(0),\ox^{22}(0)\)+\(ts_2|\ou^2|\cos\th^u,ts_2|\ou^2|\sin\th^u\)\;\mbox{ for }\;t\in[0,t_\ast),
\end{array}\right.
\end{eqnarray*}\vspace*{-0.2in}
\begin{eqnarray}\label{trajec}
\left\{\begin{array}{ll}
\ox^1(t)=\(\ox^{11}(0),\ox^{12}(0)\)+\big(ts_1|\ou^1|\cos\th^u-\eta(t-t_\ast),ts_1|\ou^1|\sin\th^u-\eta(t-t_\ast)\big)\;\mbox{ and}\\
\ox^2(t)=\(\ox^{21}(0),\ox^{22}(0)\)+\big(ts_2|\ou^2|\cos\th^u+\eta(t-t_\ast),ts_2|\ou^2|\sin\th^u+\eta(t-t_\ast)\big)\;\mbox{ for }\;t\in[t_\ast,\oT].
\end{array}\right.
\end{eqnarray}
Note that $\|\ox^2(t_\ast)-\ox^1(t_\ast)\|=R_1+R_2$ at the contacting time $t=t_{\ast}$. Recalling then that we use the sum norm yields the following equation for $t_\ast$:
\begin{eqnarray}\label{t_ast}
\begin{array}{ll}
|\ox^{21}(0)-\ox^{11}(0)+t_\ast\(s_2|\ou^2|-s_1|\ou^1|\)\cos\th^u |+
|\ox^{22}(0)-\ox^{12}(0)+t_\ast\(s_2|\ou^2|-s_1|\ou^1|\)\sin\th^u |=R_1+R_2,
\end{array}
\end{eqnarray}
which makes the connection between $t_\ast$ and the control $\ou=(\ou^1,\ou^2)$ via the given model data. Since the trajectory does not hit the boundary before $t=t_\ast$, we deduce from  {\bf 4},  {\bf 5}, and  {\bf 8} that 
\begin{equation}\label{qx}
q(t)\equiv-\mu\ox(\oT)-(1,-1)\big(\eta(\oT)+\gg([t,\oT])\big)=:(r^1_*,r^2_*)\in \R^4
\end{equation}
for all $t\in[0,t_\ast)$. Thus the maximization condition in  {\bf 6} tells us that
$$
\la r_\ast^1,\ou^1(t)\ra+\la r_\ast^2,\ou^2(t)\ra=\max_{u\in U}\big\{\la r_\ast^1,u^1\ra+\la r_\ast^2,u^2\ra\big\}
$$
for a.e.\ $t\in[0,t_\ast)$, which can be written by \eqref{qx} as
\begin{equation}\label{max_con}
\begin{split}
&\big\la -\mu\ox^1(\oT)-\mu\ox^2(\oT)+
(\eta(\oT)+\gg([t,\oT])),\big(-\ou^1(t)+\ou^2(t)\big)\big\ra\\
=&\underset{u\in U}\max\big\{\big\la -\mu\big(\ox^1(\oT)+\ox^2(\oT)\big)+
\big(\eta(\oT)+\gg([t,\oT])\big),(-u^1+u^2)\big\ra\big\}
\end{split}
\end{equation}
with $\eta(\oT)\ge 0$. Considering the case when the optimal control controls $\ou^1$ and $\ou^2$ are interior points of the control domain $U$, we have from the maximization of the function
$$
\phi(u^1,u^2):=\big\la u^1,\big(-\eta(\oT)-\gg([t,\oT])-\mu\ox^1(\oT)\big)\big\ra+\big\la u^2,\big(\eta(\oT)+\gg([t,\oT])-\mu\ox^2(\oT)\big)\big\ra
$$
the following explicit conditions:
\begin{itemize}
\item If $\|\ou^1(t)\|< 3$, then $-\eta(\oT)-\gg([t,\oT])=\mu\ox^1(\oT)$.
\item If $\|\ou^2(t)\|< 3$, then $\eta(\oT)+\gg([t,\oT])=\mu\ox^2(\oT)$.
\end{itemize}
If both cases above occur, then we get $\mu\big(\ox^1(\oT)+\ox^2(\oT)\big)=0.$ The last equation gives us the relationship
\begin{equation*}
\ox^1(\oT)+\ox^2(\oT)=0,
\end{equation*}
provided that $\mu=1$ assumed without loss of generality. Combining the latter with \eqref{trajec} gives us $\oT|\ou^2|\approx 8.3275e-13$. Furthermore, it follows from \eqref{t_ast} that either $t_\ast|\ou^2|\approx 6.1372e-13$, or $t_\ast |\ou^2|\approx 6.7832e-13$. Examine now both these possibilities for $t_\ast$.

{\bf Case 1:} If $t_\ast|\ou^2|\approx 6.1372e-13$, then we get $\oT\approx 1.3569t_\ast$. Furthermore, using \eqref{eq-eta2} and the fminsearch program from MATLAB for the cost functional $J[\ox,\ou,\oT]$, we obtain that
\begin{eqnarray*}
\eta=\frac{1}{2}\Big(2s_1|\ou^2|\Big(\frac{\sqrt{2}}{2}\Big)-s_2|\ou^2|\Big(\frac{\sqrt{2}}{2}\Big)\Big)=(1.1609e+14)|\ou^2|\ne 0
\end{eqnarray*}
and that the following expressions hold:
\begin{eqnarray*}
\left\{\begin{array}{ll}
|\ou^1|&\approx2.3806, |\ou^2|\approx1.1903,\\
\ox^1(t)&\approx\(-350-(5.3126e+14)t,-350-(5.3126e+14)t\),\quad t\in [0,5.1688e-13),\\
\ox^1(t)&\approx\(-278.5765+(3.9308e+14)t,-278.5765+(3.9308e+14)t\),\quad t\in [5.1688e-13,7.0135e-13],\\
\ox^2(t)&\approx\(-200-(2.5489e+14)t,-200-(2.5489e+14)t\),\quad t\in [0,5.1688e-13),\\
\ox^2(t)&\approx\(-271.4235+(3.9307e+14)t,-271.4235+(3.9307e+14)t\),\quad t\in[5.1688e-13,7.0135e-13],
\end{array}\right.
\end{eqnarray*}
with $t_\ast=5.1688e-13$, $\oT=7.0135e-13$, and $J=14.3596$.

{\bf Case 2:} If $t_\ast|\ou^2|\approx 6.7832e-13$, then we get $\oT\approx 1.2277t_\ast$. Furthermore, using \eqref{eq-eta2} and fminsearch program from MATLAB for the cost functional $J[\ox,\ou,\oT]$, we obtain that
\begin{eqnarray*}
\left\{\begin{array}{ll}
|\ou^1|&\approx2.3804, |\ou^2|\approx1.1902,\\
\ox^1(t)&\approx\(-350+(5.3122e+14)t,-350+(5.3122e+14)t\),\quad t\in [0,5.6744e-13),\\
\ox^1(t)&\approx\(-271.5966+(3.9305e+14)t,-271.5966+(3.9305e+14)t\),\quad t\in [5.6744e-13,6.9665e-13],\\
\ox^2(t)&\approx\(-200+(2.5487e+14)t,-200+(2.5487e+14)t\),\quad t\in [0,5.6744e-13),\\
\ox^2(t)&\approx\(-278.4034+(3.9304e+14)t,-278.4034+(3.9304e+14)t\),\quad t\in[5.6744e-13,6.9665e-13].
\end{array}\right.
\end{eqnarray*}
Then the corresponding value of the cost functional is $J=36.3565$.

Comparing the cost functional values in the two cases above, we confirm that $(\ou^1,\ou^2)=(5.3126e+14,5.3126e+14,3.9308e+14,3.9308e+14)$ is the optimal control to this problem.\vspace*{-0.05in} 

\section{Concluding Remarks}\label{sec:Conclusions}

The paper presents applications of newly derived necessary optimality conditions for free-time optimal control of sweeping processes to practical problems arising in marine surface vehicle and nanoparticles models. This is done for the first time in the literature. The obtained results and calculations demonstrate the efficiency of the established optimality conditions to solve practical problems under certain assumptions that simplify the calculations. In our future research on these and related topics, we intend to relax the imposed assumptions to be able, in particular, to determine optimal strategies in these and related models among variable control actions in the corresponding sweeping processes.

{\bf Acknowledgments.} The authors are gratefully indebted to Giovanni Colombo for his great contributions to our joint paper \cite{cmnn23} and further discussions on the material presented in this paper. We also thank Messaoud Bounkhel, a former student of Lionel Thibault, for drawing our attention to his papers \cite{bk,hb1} on the modeling and simulation in USV and nanoparticle systems.

\end{document}